\documentclass{amsproc}
\usepackage{amsmath}
\usepackage{amssymb,latexsym}
\usepackage{amscd}
\usepackage{xspace}
\usepackage{verbatim}
\begin{document}
\numberwithin{equation}{section}

\title[Initial value problems with measures]{Initial value problems for diffusion equations with singular potential}
\author{Konstantinos T. Gkikas }
\address{Centro de Modelamiento Matem\`atico, Universidad de Chile\\
Santiago de Chile, CHILE}
 \email{kgkikas@dim.uchile.cl}
\author{ Laurent Veron}
\address{Laboratoire de Math\'ematiques, Facult\'e des Sciences\\
Universit\'e Fran\c{c}ois Rabelais, Tours, FRANCE}
\email{veronl@lmpt.univ-tours.fr}

\newcommand{\norm}[1]{\left \|#1\right \|}
\newcommand{\R}{I \!  \! R}
\newcommand{\Proof}{{\it Proof. }}
\newcommand{\be}{\begin{equation}}
\newcommand{\ee}{\end{equation}}
\newcommand{\bea}{\begin{eqnarray}}
\newcommand{\eea}{\end{eqnarray}}
\def\ba{\mathbf a}
\def\bb{\mathbf b}
\newcommand{\la}{\label}
\newcommand{\up}{\upsilon}
\newcommand{\xa}{\alpha}
\newcommand{\xb}{\beta}
\newcommand{\xg}{\gamma}
\newcommand{\xG}{\Gamma}
\newcommand{\xd}{\delta}
\newcommand{\xD}{\Delta}
\newcommand{\xe}{\varepsilon}
\newcommand{\e}{\epsilon}
\newcommand{\xz}{\zeta}
\newcommand{\xh}{\eta}
\newcommand{\Th}{\Theta}
\newcommand{\xk}{\kappa}
\newcommand{\xl}{\lambda}
\newcommand{\xL}{\Lambda}
\newcommand{\xm}{\mu}
\newcommand{\xn}{\nu}
\newcommand{\ks}{\xi}
\newcommand{\KS}{\Xi}
\newcommand{\xp}{\pi}
\newcommand{\xP}{\Pi}
\newcommand{\xr}{\rho}
\newcommand{\xs}{\sigma}
\newcommand{\xS}{\Sigma}
\newcommand{\xf}{\phi}
\newcommand{\xF}{\Phi}
\newcommand{\ps}{\psi}
\newcommand{\PS}{\Psi}
\newcommand{\xo}{\omega}
\newcommand{\xO}{\Omega}
\newcommand{\C}{{\bf (C)}}
\newcommand{\CC}{ (C)}
\newcommand{\RR}{ (R)}
\newcommand{\Ren}{ I \! \! R^N}
\newcommand{\Real}{ I \! \! R}
\newcommand{\Ded}{{\cal D}^{1,2}}
\newcommand{\ho}{H_{\xO}}
\newcommand{\co}{C_{\xO}}
\newcommand{\weight}{|x|^{-(N-2)}}
\newcommand{\Lp}{L^p (\xO)}
\newcommand{\spc}{{\cal C}_0^{\infty} (\xO)}
\newcommand{\sph}{H_0^1 (\xO)}
\newcommand{\spw}{W_0^{1,2} (\xO ;\;  \weight )}
\newcommand{\iuu}{\int_{\Ren} |\nabla u|^2\,dx}
\newcommand{\ioo}{\int_{\Ren} Vu^2\,dx}
\newcommand{\Br}{B_r}
\newcommand{\bBr}{\partial \! B_r}
\newcommand{\ra}{\rightarrow}
\newcommand{\rft}{\rightarrow +\infty}
\newcommand{\bin}{\int_{\bBr}}
\newcommand{\pa}{\partial}
\newcommand{\abs}[1]{\lvert#1\rvert}
\newcommand{\darr}[4]{{\left\{\begin{array}{ll}
   {#1}&{#2}\\
   {#3}&{#4}
\end{array}\right.}}
\newcommand{\BA}{
\begin{array}}
\newcommand{\EA}{
\end{array}}
\newcommand{\ity}{\infty}
\newcommand{\prt}{\partial}
\newcommand{\opname}[1]{\mbox{\rm #1}\,}
\newcommand{\supp}{\opname{supp}}
\newcommand{\dist}{\opname{dist}}
\newcommand{\myfrac}[2]{{\displaystyle \frac{#1}{#2} }}
\newcommand{\myint}[2]{{\displaystyle \int_{#1}^{#2}}}
\newcommand{\mysum}[2]{{\displaystyle \sum_{#1}^{#2}}}
\newcommand {\dint}{{\displaystyle \int\!\!\int}}
\newcommand{\bth}[1]{\def\name{Theorem}
\begin{sub}\label{t:#1}}
\newcommand{\blemma}[1]{\def\name{Lemma}
\begin{sub}\label{l:#1}}
\newcommand{\bcor}[1]{\def\name{Corollary}
\begin{sub}\label{c:#1}}
\newcommand{\bdef}[1]{\def\name{Definition}
\begin{sub}\label{d:#1}}
\newcommand{\bprop}[1]{\def\name{Proposition}
\begin{sub}\label{p:#1}}
\newcounter{newsection}
\newtheorem{theorem}{Theorem}[section]
\newtheorem{lemma}[theorem]{Lemma}
\newtheorem{prop}[theorem]{Proposition}
\newtheorem{coro}[theorem]{Corollary}
\newtheorem{defin}[theorem]{Definition}
\newcounter{newsec} \renewcommand{\theequation}{\thesection.\arabic{equation}}
\dedicatory   {To Patrizia Pucci, with friendship and high esteem}

         \date{}
\def\ga{\alpha}     \def\gb{\beta}       \def\gg{\gamma}
\def\gc{\chi}       \def\gd{\delta}      \def\ge{\epsilon}
\def\gth{\theta}                         \def\vge{\varepsilon}
\def\gf{\phi}       \def\vgf{\varphi}    \def\gh{\eta}
\def\gi{\iota}      \def\gk{\kappa}      \def\gl{\lambda}
\def\gm{\mu}        \def\gn{\nu}         \def\gp{\pi}
\def\vgp{\varpi}    \def\gr{\rho}        \def\vgr{\varrho}
\def\gs{\sigma}     \def\vgs{\varsigma}  \def\gt{\tau}
\def\gu{\upsilon}   \def\gv{\vartheta}   \def\gw{\omega}
\def\gx{\xi}        \def\gy{\psi}        \def\gz{\zeta}
\def\Gg{\Gamma}     \def\Gd{\Delta}      \def\Gf{\Phi}
\def\Gth{\Theta}
\def\Gl{\Lambda}    \def\Gs{\Sigma}      \def\Gp{\Pi}
\def\Gw{\Omega}     \def\Gx{\Xi}         \def\Gy{\Psi}
\def\CS{{\mathcal S}}   \def\CM{{\mathcal M}}   \def\CN{{\mathcal N}}
\def\CR{{\mathcal R}}   \def\CO{{\mathcal O}}   \def\CP{{\mathcal P}}
\def\CA{{\mathcal A}}   \def\CB{{\mathcal B}}   \def\CC{{\mathcal C}}
\def\CD{{\mathcal D}}   \def\CE{{\mathcal E}}   \def\CF{{\mathcal F}}
\def\CG{{\mathcal G}}   \def\CH{{\mathcal H}}   \def\CI{{\mathcal I}}
\def\CJ{{\mathcal J}}   \def\CK{{\mathcal K}}   \def\CL{{\mathcal L}}
\def\CT{{\mathcal T}}   \def\CU{{\mathcal U}}   \def\CV{{\mathcal V}}
\def\CZ{{\mathcal Z}}   \def\CX{{\mathcal X}}   \def\CY{{\mathcal Y}}
\def\CW{{\mathcal W}} \def\CQ{{\mathcal Q}}
\def\BBA {\mathbb A}   \def\BBb {\mathbb B}    \def\BBC {\mathbb C}
\def\BBD {\mathbb D}   \def\BBE {\mathbb E}    \def\BBF {\mathbb F}
\def\BBG {\mathbb G}   \def\BBH {\mathbb H}    \def\BBI {\mathbb I}
\def\BBJ {\mathbb J}   \def\BBK {\mathbb K}    \def\BBL {\mathbb L}
\def\BBM {\mathbb M}   \def\BBN {\mathbb N}    \def\BBO {\mathbb O}
\def\BBP {\mathbb P}   \def\BBR {\mathbb R}    \def\BBS {\mathbb S}
\def\BBT {\mathbb T}   \def\BBU {\mathbb U}    \def\BBV {\mathbb V}
\def\BBW {\mathbb W}   \def\BBX {\mathbb X}    \def\BBY {\mathbb Y}
\def\BBZ {\mathbb Z}

\def\GTA {\mathfrak A}   \def\GTB {\mathfrak B}    \def\GTC {\mathfrak C}
\def\GTD {\mathfrak D}   \def\GTE {\mathfrak E}    \def\GTF {\mathfrak F}
\def\GTG {\mathfrak G}   \def\GTH {\mathfrak H}    \def\GTI {\mathfrak I}
\def\GTJ {\mathfrak J}   \def\GTK {\mathfrak K}    \def\GTL {\mathfrak L}
\def\GTM {\mathfrak M}   \def\GTN {\mathfrak N}    \def\GTO {\mathfrak O}
\def\GTP {\mathfrak P}   \def\GTR {\mathfrak R}    \def\GTS {\mathfrak S}
\def\GTT {\mathfrak T}   \def\GTU {\mathfrak U}    \def\GTV {\mathfrak V}
\def\GTW {\mathfrak W}   \def\GTX {\mathfrak X}    \def\GTY {\mathfrak Y}
\def\GTZ {\mathfrak Z}   \def\GTQ {\mathfrak Q}

\maketitle

\begin{abstract}
Let $V$ be a nonnegative locally bounded function defined in $Q_\infty:=\BBR^n\times(0,\infty)$. We study under what conditions on $V$ and on a Radon measure $\gm$ in $\mathbb{R}^d$ does it exist a function which satisfies $\partial_t u-\xD u+ Vu=0$ in $Q_\infty$ and $u(.,0)=\xm$. We prove the existence of a subcritical case in which any measure is admissible and a supercritical case where capacitary conditions are needed. We obtain a general representation theorem of positive solutions when $t V(x,t)$ is bounded and we prove the existence of an initial trace in the class of outer regular Borel measures.
\end{abstract}
\tableofcontents
\section{Introduction}
\setcounter{equation}{0}
 In this article we study the initial value problem for the heat equation
\begin{equation}\label {equ1}\begin{array} {ll}
\partial_tu-\Delta u +V(x,t)u=0\qquad&\text{in }Q_T:=\BBR^n\times(0,T)\\
\phantom{\prt_tua +V(xt)u}
u(.,0)=\gm\qquad&\text{in }\BBR^n,
\end{array}\end{equation}
where $V\in L^\infty_{loc}(Q_T)$ is a nonnegative function and $\gm$ a Radon measure in $\BBR^n$. By a (weak) solution of (\ref{equ1}) we mean a function  $u\in L_{loc}^1(\overline Q_T)$ such that $Vu\in L^1_{loc}(\overline Q_T)$, satisfying
\begin{equation}\label {equ3}\begin{array} {ll}
-\dint_{\!\!Q_T}\left(\prt_t\gf+\Gd\gf\right)udxdt+\dint_{\!\!Q_T}Vu\gf dxdt=\myint{\Gw}{}\gz d\gm
\end{array}\end{equation}
for every function $\gz\in C_c^{1,1;1}(\overline Q_T)$ which vanishes for $t= T$. Besides the singularity of the potential at $t=0$, there are two main difficulties which appear for constructing weak solutions : the growth of the measure at infinity and the concentration of the measure near some points in $\BBR^n$. Diffusion equations with singular potentials depending only on $x$ have been studied in connection with the stationnary equation (see e.g. \cite{Ka}). The particular case  of Hardy's potentials $v(x,t)=c\abs{x}^{-2}$ has been thorougly investigated  since the early work of Baras and Goldstein \cite{BG}, in connection with the problem of instantaneous blow-up. For  time dependent singular potentials most of the works are concentrated on the well posedness and the existence of a maximum principle; this is the case if $V\in L^{\infty}_tL^{\frac{n}{2},\infty}_x$, see e.g. \cite{Pie}.  In the case of time-singular potentials, a notion of non-autonomous Kato class have been introduced in \cite{SV} in order to prove that the evolution problem associated to the equation  is well posed in $L^1(\BBR^n)$. This class is the extension to diffusion equations of the Kato's class in Schr\" odinger operators. Other studies have been performed by probabilistic methods in order to analyze the $L^p-L^q$ regularizing effect \cite{Gu}.
To our knowledge, no work dealing with the initial value problems with measure data for singular operators has already been published. We present here an extension to evolution equations of a series of questions raised and solved in the case of Schr\" odinger stationary equations in  particular by \cite{An0},  \cite{An},  \cite{VY}, having in mind that one of the aim of this present work is to develop a framework adapted to the construction of the precise trace of solutions of  semilinear heat equations. This aspect will appear in a forthcoming work \cite{GV}).
\smallskip

We denote by $H(x,t)=\left(\frac{1}{4\gp t}\right)^{\frac{n}{2}}e^{-\frac{\abs{x}^2}{4t}}$ the Gaussian kernel in $\BBR^n$ and by $\BBH[\gm]$ the corresponding heat potential of a measure $\gm\in \mathfrak M(\BBR^n)$. Thus
\begin{equation}\label {equ4}\begin{array} {ll}
\BBH[\gm](x,t)=\left(\frac{1}{4\gp t}\right)^{\frac{n}{2}}\myint{}{}e^{-\frac{\abs{x-y}^2}{4t}}d\gm(y),
\end{array}\end{equation}
whenever this expression has a meaning: for example it is straightforward that if  $\gm\in  \mathfrak M(\BBR^n)$ satisfies
\begin{equation}\label {equ5}
\norm{\gm}_{\mathfrak M_{_T}}:=\myint{\BBR^n}{}e^{-\frac{\abs{y}^2}{4T}}d\abs{\gm}(y)<\infty,
\end{equation}
then (\ref{equ4}) has a meaning as long as $t< T$, and let be $\mathfrak M_{_T}(\BBR^n)$ the set of Radon measures in $\BBR^n$ satisfying (\ref{equ5}). If $G\subset\BBR^n$, let $Q^G_T$ be the cylinder $G\times (0,T)$, $B_R(x)$ the ball of center $x$ and radius $R$ and $B_R=B_R(0)$. We prove \smallskip

\noindent {\bf Theorem A }{\it Let the  measure $\gm$ verifies
\begin{equation}\label {equ6}\begin{array} {ll}
\dint_{\!\!Q^{B_R}_T}\BBH[\abs{\gm}](x,t)V(x,t)dxdt\leq M_{R}\qquad\forall R>0.
\end{array}\end{equation}
Then (\ref{equ1}) admits a solution in $Q_T$.}\smallskip

A measure which satisfies (\ref{equ6}) is called an {\it admissible measure} and a measure for which there exists a solution  to problem (\ref{equ1}) is called a {\it good measure}. Notice that even when $V=0$, uniqueness without any restriction on $u$ is not true, however the next uniqueness result holds:\smallskip

\noindent {\bf Theorem B }{\it Let $u$ be a weak solution of (\ref{equ1}) with $\gm=0$. If $u$ satisfies
\begin{equation}\label {equ7}\begin{array} {ll}
\dint_{\!\!Q_T}\left(1+V(x,t)\right)e^{-\gl\abs{x}^2}\abs {u(x,t)}dxdt<\infty
\end{array}\end{equation}
for some $\gl>0$, then $u=0$}. \smallskip

We denote by $\CE_\gn(Q_T)$ the set of functions $u\in L^1_{loc}(Q_T)$ for which (\ref{equ7}) holds for some $\gl>0$. The general result we prove is the following.\smallskip

\noindent {\bf Theorem C }{\it Let $\gm\in \mathfrak M(\BBR^n)$ be an admissible measure satisfying (\ref{equ5}). Then there exists a unique solution $u_\gm\in \CE_\gn(Q_T)$ to problem  (\ref{equ1}). Furthermore
\begin{equation}\label {equ8}\begin{array} {ll}
\dint_{\!\!Q_T}\left(\myfrac{n}{2T}+V\right)\abs{u}e^{-\frac{\abs{x}^2}{4(T-t)}}dxdt\leq
\myint{\BBR^n}{}e^{-\frac{\abs{y}^2}{4T}}d\abs{\gm}(y).
\end{array}\end{equation}
}\medskip

We consider first the {\it subcritical case}, which means that  any positive measure satisfying (\ref{equ5}) is a good measure and we prove that such is the case if for any $R>0$ there exist $m_R>0$ such that
\begin{equation}\label {equ9}\begin{array} {ll}
\dint_{\!\!Q^{B_R}_T}H(x-y,t)V(x,t)dxdt\leq m_Re^{-\frac{\abs{y}^2}{4T}}.
\end{array}\end{equation}
Moreover we prove a stability result among the measures satisfying (\ref{equ5}): if $V$ verifies for all $R>0$
\begin{equation}\label {equ10}\begin{array} {ll}
\displaystyle\sup_{y\in\BBR^n}e^{\frac{\abs{y}^2}{4T}}\dint_{\!\!E}H(x-y,t)V(x,t)dxdt\to 0\quad\text {when }\abs{E}\to 0\,,\, E\text { Borel subset of }Q^{B_R}_T,
\end{array}\end{equation}
then if $\{\gm_k\}$ is a sequence of Radon measures bounded in $\mathfrak M_{_T}(\BBR^n)$ which converges in the weak sense of measures to $\gm$, then $\left\{\left(u_{\gm_k},Vu_{\gm_k}\right)\right\}$ converges to $(u_{\gm},Vu_{\gm})$
 in $L^1_{loc} (\overline Q_T)$.\smallskip

 In the {\it supercritical case}, that is when not all measure in $\mathfrak M_{_T}(\BBR^n)$ is a good measure, we develop a capacitary framework in order to characterize the good measures. We denote by $\mathfrak M^V(\BBR^n)$ the set of Radon measures such that $V\BBH[\gm]\in L^1(Q_T)$ and $ \norm\gm_{\mathfrak M^V}:=\norm{V\BBH[\gm]}_{L^1}$. If $E\subset Q_T$ is a Borel set, we set
 \begin{equation}\label {equ11}\begin{array} {ll}
C_V(E)=\sup\{\gm(E):\gm\in \mathfrak M_+^V(\BBR^n),\gm(E^c)=0,\norm\gm_{\mathfrak M^V}\leq 1\}.
\end{array}\end{equation}
This defines a capacity. If
 \begin{equation}\label {equ11*}\begin{array} {ll}
C^*_V(E)=\inf\{\norm f_{L^\infty}: \check H[f](y)\geq 1\quad \forall y\in E\},
\end{array}\end{equation}
where
 \begin{equation}\label {equ12}\begin{array} {ll}
\check H[f](y)=\dint_{\!\!Q_T}H(x-y,t)V(x,t)f(x,t) dxdt=\myint{0}{T}\BBH[Vf](y,t)dt\qquad\forall y\in\BBR^n,
\end{array}\end{equation}
then $C^*_V(E)=C_V(E)$ for any compact set. Denote by $Z_V$ {\it the  singular set of $V$}, that is the largest set with zero
$C_V$ capacity. Then
 \begin{equation}\label {equ13}\begin{array} {ll}
Z_V=\{x\in\BBR^n:\dint_{\!\!Q_T}H(x-y,t)V(y,t) dxdt=\infty\},
\end{array}\end{equation}
and the following result characterizes the good measures.\smallskip

\noindent {\bf Theorem D }{\it If $\gm$ is an admissible measure then $\gm(Z_V)=0$. If $\gm\in \mathfrak M_{_T}(\BBR^n)$ satisfies $\gm(Z_V)=0$, then it is a good measure. Furthermore $\gm$ is a positive good measure if and only if there exists an increasing sequence of positive admissible measures $\{\gm_k\}$ which converges to $\gm$ in the weak $* $ topology.
}\medskip

 Since many important applications deal with the nonlinear equation
\begin{equation}\label {equ14}\begin{array} {ll}
\partial_tu-\Delta u +\abs{u}^{q-1}u=0\qquad&\text{in }Q_\infty:=\BBR^n\times(0,\infty),
\end{array}\end{equation}
where $q>1$ and due to the fact that   {\it any} solution defined in $Q_\infty$ satisfies
\begin{equation}\label {equ15}\begin{array} {ll}
\abs{u(x,t)}^{q-1}\leq \myfrac{1}{t(q-1)}\qquad\forall (x,t)\in Q_\infty,
\end{array}\end{equation}
we shall concentrate on potentials $V$ which satisfy
\begin{equation}\label {equ16}\begin{array} {ll}
0\leq V(x,t)\leq \myfrac{C_1}{t}\qquad\forall (x,t)\in Q_T,
\end{array}\end{equation}
for some $C_1>0$. For such potentials we prove the existence of a representation theorem for positive solutions of
\begin{equation}\label {equ17}\begin{array} {ll}
\partial_tu-\Delta u +V(x,t)u=0\qquad&\text{in }Q_T.
\end{array}\end{equation}

If $u$ is a positive solution of (\ref{equ1}) in $Q_T$ with $\gm\in \mathfrak M_+(\BBR^n)$, it is the increasing limit of the solutions $u=u_R$ of
\begin{equation}\label {equ18}\begin{array} {ll}
\partial_tu-\Delta u +V(x,t)u=0\qquad&\text{in }Q^{B_R}_T\\
\phantom{\partial_tu-\Delta  +V(x,t)u}u=0\qquad&\text{in }\prt B_R\times (0,T)\\
\phantom{\partial_tu-\Delta  +Vtu}u(.,0)=\chi_{B_R}\gm\qquad&\text{in }B_R,
\end{array}\end{equation}
when $R\to\infty$, thus there exists a positive function $H_V\in C(\BBR^n\times\BBR^n\times (0,T))$ such that
\begin{equation}\label {equ19}\begin{array} {ll}
u(x,t)=\myint{\BBR^N}{}H_V(x,y,t)d\gm (y).
\end{array}\end{equation}
Furthermore we show how to construct $H_V$ from $V$ and we prove the following formula
\begin{equation}\label {equ19-1}\begin{array} {ll}
H_V(x,y,t)=\myint{\BBR^N}{}e^{\psi (x,t)}\Gamma (x,\xi,t)d\xm_y(\xi),
\end{array}\end{equation}
where $\xm_y$ is a Radon measure such that
\be
\xd_y\geq \xm_y,\label{111}
\ee
($\xd_y$ is the Dirac measure concentrated at $y$),
\begin{equation}\label {equ19-2}\begin{array} {ll}
\psi(x,t)=\myint{t}{T}\myint{\BBR^n}{}\left(\myfrac{1}{4\gp (s-t)}\right)^{\frac{n}{2}}e^{-\frac{\abs{x-y}^2}{4 (s-t)}}V(y,s) dy ds
\end{array}\end{equation}
and $\Gamma$ satisfies the following estimate
\begin{equation}\label {equ19-3}\begin{array} {ll}
c_1t^{-\frac{n}{2}}e^{-\gamma_1\frac{\abs{x-y}^2}{t}}\leq \Gamma (x,y,t)\leq c_2t^{-\frac{n}{2}}e^{-\gamma_2\frac{\abs{x-y}^2}{t}}
\end{array}\end{equation}
where $A_i$, $c_i$ depends on $T$, $d$ and $V$. Conversely, we first prove the following representation result\medskip

\noindent {\bf Theorem E }{\it Assume $V$ satisfies (\ref{equ16}). If $u$ is a positive solution of (\ref{equ1}) in $Q_T$, there exists a positive Radon measure $\gm$ in $\BBR^n$ such that  (\ref{equ19}) holds.}\medskip

If $\gm\in\mathfrak M_{_T}(\BBR^n)$ is positive, we can define for any $k>0$ the solution $u_k$ of
\begin{equation}\label {equ20_*}\begin{array} {ll}
\partial_tu-\Delta u +V_k(x,t)u=0\qquad&\text{in }Q_T\\
\phantom{\prt_tua +\!V(x,t)u}
u(.,0)=\gm\qquad&\text{in }\BBR^n,
\end{array}\end{equation}
where $V_k(x,t)=\min\{k,V(x,t)\}$, and
\begin{equation}\label {equ21_*}\begin{array} {ll}
u_k(x,t)=\myint{\BBR^N}{}H_{V_k}(x,y,t)d\gm (y).
\end{array}\end{equation}
Moreover $\{H_{V_k}\}$ and $\{v_k\}$ decrease respectively to $H_{V}$ and $u^*$ there holds
\begin{equation}\label {equ22}\begin{array} {ll}
u^*(x,t)=\myint{\BBR^N}{}H_V(x,y,t)d\gm (y).
\end{array}\end{equation}
 However $u^*$ is not a solution of (\ref{equ1}), but of a relaxed problem where $\gm$ is replaced by a smaller measure $\gm^*$  called {\it  the reduced measure associated to }$\gm$.  If we define the {\it zero set of $V$} by

\begin{equation}\label {equ20}\begin{array} {ll}
\CS ing_V:=\{y\in \BBR^N:H_V(x,y,t)=0\},
\end{array}\end{equation}
we prove  \smallskip

\noindent {\bf Theorem F }{\it If

\be
\limsup_{t\rightarrow0}\myint{t}{T}\myint{\BBR^n}{}\left(\myfrac{1}{4\gp (s-t)}\right)^{\frac{n}{2}}e^{-\frac{\abs{\xi-y}^2}{4 (s-t)}}V(y,s) dy ds=\infty,\label{222*}
\ee
then $$\xi\in\CS ing_V,\; \mathrm{i.e.} \;H_V(x,\xi,t)=0,\;\forall\; (x,t)\in\mathbb{R}^n\times(0,\infty). $$}
\noindent We note here that if $V$ satisfies (\ref{222*}) then $\xd_\xi$ is not admissible measure and the reduced measure $(\xd_\xi)^*=\xm_\xi$  associated to $\xd_\xi$ is zero.  \smallskip

\noindent {\bf Theorem G }{\it Assume $V$ satisfies (\ref{equ15}) and $\gm\in\mathfrak M_{_T}(\BBR^n)$. Then\smallskip

\noindent (i) $supp (\gm-\gm^*)\subset \CS ing_V$.\smallskip

\noindent (ii) If $\gm(\CS ing_V)=0$, then $\gm^*=0$.\smallskip

\noindent (iii) $\CS ing_V= Z_V.$
}\medskip

The last section is devoted to the initial trace problem: to any positive solution $u$ of (\ref{equ1}) we can associate an open subset $\CR(u)\subset\BBR^n$ which is the set of points $y$ which possesses a neighborhood $U$ such that
\begin{equation}\label {equ21}\begin{array} {ll}
\dint_{\!\!Q^U_T}V(x,t)u(x,t)dxdt<\infty.
\end{array}\end{equation}
There exists a positive Radon measure $\gm_u$ on $\CR(u)$ such that
\begin{equation}\label {equ22^*}\begin{array} {ll}
\lim_{t\to 0}\myint{\BBR^n}{}u(x,t)\gz(x) dx=\myint{\BBR^n}{}\gz d\gm\qquad\forall \gz\in C_c(\CR(u)).
\end{array}\end{equation}
The set $\CS(u)=\BBR^n\setminus \CR(u)$ is the set of points $y$ such that for any open set $U$ containing $y$, there holds
\begin{equation}\label {equ23}\begin{array} {ll}
\dint_{\!\!Q^U_T}V(x,t)u(x,t)dxdt=\infty.
\end{array}\end{equation}
If $V$ satisfies (\ref{equ17}), $\CS(u)$ it has the property that
\begin{equation}\label {equ24}\begin{array} {ll}
\limsup_{t\to 0}\myint{U}{}u(x,t) dx=\infty.
\end{array}\end{equation}
Furthermore, if is satisfies (\ref{equ10}), then $\CS(u)=\emptyset$.\smallskip

An alternative construction of the initial trace based on the sweeping method is also developed.\smallskip

Precise definitions of the different notions used in the introduction will be given in the next sections.
\smallskip

\noindent{\bf Aknowledgements} This work has been prepared while the first author was visiting the Laboratory of Mathematics and Theoretical Physics, CNRS-UMR 7350, thanks to the support of a grant from R\'egion Centre, in the framework of the program {\it  Cr\'eation et propagations de singularit\'es dans les \'equations non-lin\'eaires}.

\section{The subcritical case}
\setcounter{equation}{0}
Let $Q_T=\mathbb{R}^n\times(0,T]$. In this section we  consider the linear parabolic problem
\begin{equation}\label{sub1}\BA {lll}
\prt_tu-\Gd u+Vu=0\qquad &\text{in }\;\; Q_T\\
\phantom{\prt_tu-\Gd ii}
u(.,0)=\xm\qquad &\text{in }\;\; \mathbb{R}^n\times\{0\},
\EA\end{equation}
where $V\in L^1_{loc}(Q_T)$ is nonnegative and $\xm$ is a Radon measure.
\begin{defin}
We say that $\xm\in\mathfrak M(\mathbb{R}^n)$ is a good measure if problem (\ref{sub1}) has a weak solution $u$ i.e.
there exists a function $u\in L^1_{loc}(\overline Q_T),$ such that $Vu\in L^1_{loc}(\overline Q_T)$ which satisfies
\be
-\int\!\!\int_{Q_T} u(\prt_t\xf+\xD\xf)dxdt+\int\!\!\int_{Q_T} Vu\xf dxdt=\int_{\mathbb{R}^n}\xf(x,0)d\xm\quad\forall\xf\in X(Q_T),\label{weakun}
\ee
where $X(Q_T)$ is the space of test functions defined by
$$X(Q_T)=\{\xf\in C_c(\overline{Q}_T),\;\prt_t\xf+\xD\xf\in L^\infty_{loc}(\overline{Q}_\infty),\;\xf(x,T)=0\}$$
\end{defin}
\begin{defin}
Let $H(x,t)$ be the heat kernel of heat equation in $\mathbb{R}^n$,  we say that $\xm\in\mathfrak M(\mathbb{R}^n)$ is an admissible measure if \\
(i) $$||V\mathbb{H}[|\xm|]||_{L^1(Q_T^{B_R})}=\int\!\!\int_{Q_T^{B_{R}}}\left(\int_{\mathbb{R}^n} H(x-y,t)d|\xm(y)|\right)V(x,t)dxdt<M_{R,T}$$
where $M_{R,T}$ is a positive constant.
\end{defin}

\begin{defin}
A function $u(x,t)$ will be said to belong to the class $\mathcal{E_V}(Q_T)$ if there exists $\xl>0$ such that
$$\int\!\!\int_{Q_T}e^{-\xl|x|^2}|u(x,t)|(1+V(x,t))dxdt<\infty.$$
A measure in $\BBR^n$ belongs to the class $\mathfrak M_{_T}(\BBR^n)$ if
 $$\norm{\gm}_{\mathfrak M_T}:=\int_{\mathbb{R}^n}e^{-\frac{|x|^2}{4T}}d|\xm|<\infty.$$
\end{defin}

\begin{lemma}
There exists at most one weak solution of problem (\ref{sub1}) in the class $\mathcal{E_V}(Q_T).$\label{monad}
\end{lemma}
\Proof
Let  $u_1$ and $u_2$ be two solutions in the class $\mathcal{E_V}(Q_T)$ then $w=u_1-u_2$ is a solution with initial data 0. Choose a standard mollifier $\xr:B(0,1)\mapsto [0,1]$ and define
$$w_j(x,t)=j^n\int_{B_{\frac{1}{j}}(x)}\xr(j(x-y))w(y,t)dy\equiv\int_{B_{\frac{1}{j}}(x)}\xr_j(x-y)w(y,t)dy.$$
Then $w_j(.,t)$ is $C^\infty$ and from the equation satisfied by $w$, it holds
$$\prt_tw_j-\xD w_j+\int_{B_{\frac{1}{j}}(x)}V(y,t)\xr_j(x-y)w(y)dy=0,$$
where $\prt_tw_j$ is taken in the weak sense. \smallskip

\noindent First we consider the case $\xl>0$ and $t\leq\min\{\frac{1}{16\xl},T\}=T'.$\\
Set $\xf(x,t)=\xi(x,t)\xz(x),$ where $\xi(x,t)=e^{-\frac{|x|^2}{4(\frac{1}{8\xl}-t)}}$ and $\xz\in C_c^\infty(\mathbb{R}^n).$
Given $\xe>0$ we define
$$g_j=\sqrt{w_j^2+\xe}.$$
Because
$\prt_t(g_j\xf)=\frac{\prt_tw_j}{\sqrt{w_j^2+\xe}}\xf+g_j\prt_t\xf,$, by a straightforward calculation we have
\bea
\nonumber
\int_{\mathbb{R}^n}\left[g_j\xf(.,s)\phantom{\frac{x}{x}}\!\!\!\!\!\right]_{s=t}^{s=0} dx&=&\int\!\!\int_{Q_t}\frac{w_j}{\sqrt{w_j^2+\xe}}\xf\xD w_jdxds\\ \nonumber
&-&\int\!\!\int_{Q_t}\frac{w_j(x,s)}{\sqrt{w_j^2(x,s)+\xe}}\xf(x,s)
\left(\int_{B_{\frac{1}{j}}(x)}V(y,t)\xr_j(x-y)w(y,s)dy\right)dxds\\ \nonumber
&+&\int\!\!\int_{Q_t}g_j\xf_sdxds\\ \nonumber
&=& I_1+I_2+I_3.
\eea
By integration by parts, we obtain
\bea
\nonumber
I_1&=&-\int\!\!\int_{Q_t}\frac{|\nabla w_j|^2}{\sqrt{w_j^2(x,s)+\xe}}\xf dxds+\int\!\!\int_{Q_t}\frac{|\nabla w_j|^2 w_j^2}{(w_j^2(x,s)+\xe)^\frac{3}{2}}\xf dxds-\int\!\!\int_{Q_t}\frac{w_j}{\sqrt{w_j^2+\xe}}\nabla w_j.\nabla\xf dxds\\ \nonumber
&\leq&-\int\!\!\int_{Q_t}\frac{w_j}{\sqrt{w_j^2+\xe}}\nabla w_j.\nabla\xf dxds\\ \nonumber
&\leq&-\int\!\!\int_{Q_t}\nabla g_j.\nabla\xf dxds\\ \nonumber
&=&-\int\!\!\int_{Q_t}\xz\nabla g_j.\nabla\xi dxds-\int\!\!\int_{Q_t}\xi\nabla g_j.\nabla\xz dxds\\ \nonumber
&=&\int\!\!\int_{Q_t}\xz g_j\xD\xi dxds+\int\!\!\int_{Q_t}g_j\nabla\xz.\nabla\xi dxds.
\eea
Since $t\leq T$, there holds $\xi|\nabla g_j|\in L^1(Q_{T'}),$  $\xi g_j\in L^1(Q_{T'}),$ $|\xD\xi| g_j\in L^1(Q_{T'}),$ $\prt_{s}\xi g_j\in L^1(Q_{T'})$ and
$$\int\!\!\int_{Q_t}\frac{w_j(x,s)}{\sqrt{w_j^2(x,s)+\xe}}
\left(\int_{B_{\frac{1}{j}}(x)}V(y,t)\xr_j(x-y)w(y,s)dy\right)\xi dxds<\infty.$$
The reason for which $\xi|\nabla g_j|\in L^1(Q_{T'})$  follows from the next inequality
$$\BA {l}\displaystyle\int\!\!\int_{Q_{T'}}|\nabla g_j|\xi dxds=\int\!\!\int_{Q_{T'}}\frac{|\nabla w_j|}{\sqrt{\xe+w^2_j}}\xi dxds\\[4mm]
\phantom{\displaystyle\int\!\!\int_{Q_{T'}}|\nabla g_j|\xi dxds}
\leq \displaystyle\int\!\!\int_{Q_{T'}}e^{-\frac{|x|^2}{4(\frac{1}{8\xl}-t)}}\left(\int_{B_{\frac{1}{j}}(x)}|\nabla\xr_j(x-y)|w(y,s)dy\right)dxds.\EA$$
Since $\forall y\in B_{\frac{1}{j}}(x)$, we have $|x|^2\geq(|y|^2-\frac{1}{j})^2=|y|^2+\frac{1}{j^2}-2\frac{|y|}{j}\geq\frac{|y|^2}{2}-(C-1)\frac{1}{j^2},$ for some positive constant $C>0$ independent on $j,y$ and $x$. Thus we have, using the fact that $e^{-\xl|y|^2}w\in L^1(Q_{T})$,
$$\int\!\!\int_{Q_{T'}}|\nabla g_j|\xi dxds\leq C(j,\xl)\int\!\!\int_{Q_{T'}}\int_{B_{\frac{1}{j}}(x)}e^{-\frac{|y|^2}{8(\frac{1}{8\xl}-t)}}
|\nabla\xr_j(x-y)|w(y,s)dydxds<\infty.$$
Also
$$
\int\!\!\int_{Q_t}\frac{w_j(x,s)}{\sqrt{w_j^2(x,s)+\xe}}\xi
\left(\int_{B_{\frac{1}{j}}(x)}V(y,t)\xr_j(x-y)w(y,s)dy\right)dxds\to_{j\to\infty}
\int\!\!\int_{Q_t}\frac{w^2(x,s)}{\sqrt{w^2(x,s)+\xe}}\xi V(y,t)dxds
$$
and
$$\int_{\mathbb{R}^n}\sqrt{w_j^2(x,s)+\xe}(\xi_s+\xD\xi) dxds\to_{j\to\infty}
\int_{\mathbb{R}^n}\sqrt{w^2(x,s)+\xe}(\xi_s+\xD\xi) dxds.
$$
We choose $\xz_R=1$ in $B_R,$ $0\leq\xz_R\leq1$ in $B_{R+1}\setminus B_R$ and 0 otherwise. Letting successively $j\to\infty$, $R\to\infty$ and finally $\xe\to0$, we derive
$$\int_{\mathbb{R}^n}|w(x,t)|\xi(x,t)dx\leq\int\!\!\int_{Q_t}|w|(\xi_s+\xD\xi) dxds-\int\!\!\int_{Q_t}w(x,s)\xi V(y,t)dxds.$$
Since
 $$\xi_s+\xD\xi=-\frac{n}{2(\frac{1}{8\xl}-s)},$$
 and $V\geq0,$ we have $w(x,t)=0\;\forall\;(x,t)\in Q_{{T'}}$. If $T'=T$ this complete the proof for $\xl\geq0$, otherwise the proof can be completed by a finite number of interations of the same argument on $\mathrm{R}^n\times(T',2T'),$ $\mathrm{R}^n\times(2T',3T')$, etc. If $\xl=0$ we set $\xi=1$ and the result follows by similar argument\hfill$\Box$

\begin{theorem}\label{fragma}
If $\xm\in \mathfrak M_{_T}(\BBR^n)$ is an admissible measure, there exists a unique $u=u_\gm\in \mathcal{E_V}(Q_T)$  solution of
(\ref{sub1}). Furthermore the following estimate holds
\begin{equation}\label{intest}
\frac{n}{2T}\int\!\!\int_{Q_T}|u|e^{-\frac{|x|^2}{4(T-t)}}dxds
+\int\!\!\int_{Q_T}|u|Ve^{-\frac{|x|^2}{4(T-t)}}dxds
\leq\int_{\mathbb{R}^n}e^{-\frac{|x|^2}{4T}}d|\xm|.
\end{equation}
\end{theorem}
\Proof
First we assume that $\xm\geq0.$ Let $\xm_R=\chi_{B_R}\xm.$ It is well known that the heat kernel $H^{B_R}(x,y,t)$  in $\xO=B_R$ is increasing with respect to $R$ and $H^{B_R}\to H,$ as $R\to\infty$ in $L^1(Q_T)$ for any $T>0$. Thus $\xm_R$ is an admissible measure in $B_R$ and by Proposition \ref{subcrit}, there exists a unique weak solution $u_R$ of problem \ref{W2} on $\xO=B_R.$ By (ii) of Proposition \ref{prop} we have
$$
-\int\!\!\int_{Q_T}|u_R|(\prt_t\xf+\xD\xf)dxdt+\int\!\!\int_{Q_T} |u_R|V\xf dxdt\leq\int_{B_R}\xf(x,0)d|\xm_R|.
$$
If we set $\xf_\xe(x,t)=e^{-\frac{|x|^2}{4(T+\xe-t)}};\;\xe>0,$ then
$$\prt_t\xf+\xD\xf=-\frac{n}{2(T+\xe-t)}e^{-\frac{|x|^2}{4(T+\xe-t)}},$$
thus we have
$$\int\!\!\int_{Q_T}|u_R|\frac{n}{2(T+\xe-t)}e^{-\frac{|x|^2}{4(T+\xe-t)}}dxdt+\int\!\!\int_{Q_T} |u_R|Ve^{-\frac{|x|^2}{4(T+\xe-t)}} dxdt\leq \int_{B_R}e^\frac{-|x|^2}{4T+4\xe}d\xm_R,$$
which implies
$$\frac{n}{2T+\xe}\int_0^T\int_{B_R}|u_R|e^{-\frac{|x|^2}{4(T+\xe-t)}}dxdt+\int_0^T\int_{B_R} |u_R|Ve^{-\frac{|x|^2}{4(T+\xe-t)}} dxdt\leq\int_{\mathbb{R}^n}e^\frac{-|x|^2}{4T+4\xe}d\xm_R.$$
Letting $\xe\to0,$ we derive
$$\frac{n}{2T}\int\!\!\int_{Q_T}|u_R|e^{-\frac{|x|^2}{4(T-t)}}dxdt+\int\!\!\int_{Q_T} |u_R|Ve^{-\frac{|x|^2}{4(T-t)}} dxdt\leq\int_{\mathbb{R}^n}e^\frac{-|x|^2}{4T}d\xm_R\leq\int_{\mathbb{R}^n}e^\frac{-|x|^2}{4T}d\xm.$$
Now by the maximum principle $\{u_R\}$ is increasing with respect to $R$  and converges to some function $u.$ By the above inequality $u\in\mathcal{E_V}(Q_T)$ satisfies the estimate (\ref{fragma}) and $u$ is a weak solution of problem (\ref{sub1}). By Lemma \ref{monad} it is unique. In the general case we write $\xm=\xm^+-\xm^-$ and the result follows by the above arguments and Lemma \ref{monad}. {\it In the sequel we shall denote by $u_\gm$ this unique solution}.\hfill$\Box$\medskip

\begin{defin}
A potential $V$ is called subcritical in $Q_T$ if for any $R>0$ there exists $m_R>0$ such that
\begin{equation}\label{subcritV}
\dint_{\!\!Q_T^{B_R}}H(x-y,t)V(x,t)dxdt\leq m_Re^{-\frac{\abs{y}^2}{4T}}\qquad\forall y\in\BBR^n.
\end{equation}
It is called strongly subcritical if moreover
\begin{equation}\label{stabV}
e^{\frac{\abs{y}^2}{4T}}\dint_{\!\!E}H(x-y,t)V(x,t)dxdt\to 0\quad\text {when }\abs{E}\to 0\,,\, E\text { Borel subset of }Q^{B_R}_T,
\end{equation}
uniformly with respect to $y\in\BBR^n$
\end{defin}

 \begin{theorem}\label{uncondi} Assume $V$ is subcritical. Then any measure in $\mathfrak M_T(\BBR^n)$ is admissible. Furthermore, if $V$ is strongly subcritical and $\{\gm_k\}$ is a sequence of measures uniformly bounded in $\mathfrak M_T(\BBR^N)$ which converges weakly to $\gm$, then the corresponding solutions  $\{u_{\gm_k}\}$ converge to $u_{\gm}$ in $L^1_{loc}(\overline Q_T)$, and  $\{Vu_{\gm_k}\}$ converges to $Vu_\gm$ in $L^1_{loc}(\overline Q_T)$.
 \end{theorem}
\Proof For the first statement we can assume $\gm\geq 0$ and there holds
$$\BA{ll}\displaystyle\dint_{\!\!Q_T^{B_R}} H(x-y,t)d\xm(y)V(x,t)dxdt=
\int_{\BBR^n}\left(\dint_{\!\!Q_T^{B_R}} H(t,x-y)V(x,t)dxdt\right) d\xm(y)\\[4mm]
\phantom{\displaystyle\dint_{\!\!Q_T^{B_R}} H(x-y,t)d\xm(y)V(x,t)dxdt}
\leq \displaystyle m_R\int_{\BBR^n}e^{-\frac{\abs{y}^2}{4T}}d\gm(y)\\[4mm]
\phantom{\displaystyle\dint_{\!\!Q_T^{B_R}} H(x-y,t)d\xm(y)V(x,t)dxdt}
\leq m_R\norm{\gm}_{\mathfrak M_T}.
\EA$$
Thus $\gm$ is admissible. For the second  statement, we assume first that $\gm_k\geq 0$. By lower semicontinuity $\gm\in \mathfrak M_T(\BBR^N)$ and $\norm{V\BBH[\gm]}_{L^1(Q_T^{B_R})}\leq M_{R,T}$ for any $k$. Since $0\leq u_{\gm_k}\leq \BBH[\gm_k]$ and $\BBH[\gm_k]\to \BBH[\gm]$ in $L^1_{loc}(\overline Q_T)$, the sequence $\{u_{\gm_k}\}$ is uniformly integrable and thus relatively compact in $L^1_{loc}(\overline Q_T)$. Furthermore  $0\leq Vu_{\gm_k}\leq V\BBH[\gm_k]$. Let $E\subset Q_T^{B_R}$ be a Borel subset, then
$$\BA {l}\displaystyle\int\!\!\int_EV\BBH[\gm_k]dxdt=\int_{\BBR^n}\left(\int\!\!\int_EVH(x-y,t)dxdt\right)d\gm_k(y)\\[4mm]
\phantom{\displaystyle\int\!\!\int_EV\BBH[\gm_k]dxdt}=\displaystyle\int_{\BBR^n}\left(e^{\frac{\abs{y}^2}{4T}}\int\!\!\int_EV(x)H(x-y,t)dxdt\right)e^{-\frac{\abs{y}^2}{4T}}d\gm_k(y)\\[4mm]
\phantom{\displaystyle\int\!\!\int_EV\BBH[\gm_k]dxdt}\leq \displaystyle\ge(\abs{E})\norm{\gm_k}_{\mathfrak M_T},
\EA$$
where $\ge(r)\to 0$ as $r\to 0$. Thus $\{(u_{\gm_k},Vu_{\gm_k})\}$ is locally compact in $L^1_{loc}(\overline Q_T)$ and, using a diagonal sequence, there exist $u\in L^1_{loc}(\overline Q_T)$ with $Vu\in L^1_{loc}(\overline Q_T)$ and a subsequence $\{k_j\}$ such that
$\{(u_{\gm_{k_j}},Vu_{\gm_{k_j}})\}$ converges to $(u_\gm, Vu_\gm)$ a.e. and in $L^1_{loc}(\overline Q_T)$. From the integral expression (\ref{weakun}) satisfied by the $u_{\gm_k}$, $u$ is a weak solution of problem (\ref{sub1}). Since the $u_{\gm_k}$ satisfy (\ref{intest}), the property holds for $u$, therefore$u=u_\gm$ is the unique solution of (\ref{sub1}), which ends the proof.\hfill$\Box$.\medskip

As a variant of the above result which will be useful later on we have

\begin{prop}\label{uncondi+} Assume $V$ satisfies
\begin{equation}\label{stabV+}
e^{\frac{\abs{y}^2}{4T}}\dint_{\!\!E}H(x-y,t)V(x,t+\gt)dxdt\to 0\quad\text {when }\abs{E}\to 0\,,\, E\text { Borel subset of }Q^{B_R}_T,
\end{equation}
uniformly with respect to $y\in\BBR^n$ and $\gt\in [0,\gt_0]$. Let $\gt_k>0$ with $\gt_k\to 0$ and $\{\gm_k\}$ be a sequence uniformly bounded in $\mathfrak M_T(\BBR^N)$ which converges weakly to $\gm$. Then the solutions  $\{u_{\gt_k,\gm_k}\}$ of
\begin{equation}\label{equV+}\BA {ll}
\prt_tu-\Gd u+Vu=0\quad&\text{on }\;\;\BBR^n\times (\gt_k,T)\\\phantom{-\Gd u+.}
u(.,\gt_k)=\gm_k&\text{on }\;\;\BBR^n\times \{\gt_k\}
\EA\end{equation}
(extended by $0$ on $(0,\gt_k)$) converge to $u_{\gm}$ in $L^1_{loc}(Q_T)$, and  $\{Vu_{\gm_k}\}$ converges to $Vu_\gm$ in $L^1_{loc}( Q_T)$.
\end{prop}

Condition (\ref{stabV}) may be very difficult to verify and we give below a  sufficient condition for it to hold.

\begin{prop}\label{suff}
Assume $V$ satisfies
\begin{equation}\label{stabV'}
\lim_{\gl\to 0}e^{\frac{\abs{y}^2}{4T}}\gl^{-n}\myint{0}{\gl}\myint{B_{\gl^2}(y)}{}V(x,t)dxdt=0
\end{equation}
uniformly with respect to $y\in\BBR^n$, then $V$ is strongly subcritical.
\end{prop}
\Proof Let $E\subset Q_T^{B_R}$ be a Borel set. For $\gd>0$, we define the weighted heat ball of amplitude $\gd e^{-\frac{\abs{y}^2}{4T}}$ by
$$P_\gd=P_\gd(y,T)=\left\{(x,t)\in Q_T:H(x-y,t)\geq\gd e^{-\frac{\abs{y}^2}{4T}}\right\}.$$
By an straightforward computation, one sees that
$$P_\gd(y,T)\subset B_{a_n\gd^{-\frac{1}{n}}e^{\frac{\abs{y}^2}{4nT}}}(y)\times[0,b_n\gd^{-\frac{2}{n}}e^{\frac{\abs{y}^2}{2nT}}]:=R_\gd(y,T),
$$
for some $a_n,b_n>0$. We write
$$
\dint_{\!\!E}H(x-y,t)V(x,t)dxdt=\dint_{\!\!E\cap P_\gd}H(x-y,t)V(x,t)dxdt+\dint_{\!\!E\cap P^c_\gd}H(x-y,t)V(x,t)dxdt.
$$
Then
$$\dint_{\!\!E\cap P^c_\gd}H(x-y,t)V(x,t)dxdt\leq \gd e^{-\frac{\abs{y}^2}{4T}}\dint_{\!\!E}V(x,t)dxdt,
$$
and
$$\BA {l}\dint_{\!\!E\cap P_\gd}H(x-y,t)V(x,t)dxdt\leq \myint{0}{\gd}\myint{\{(x,t)\in Q_T^{B_R}:H(x-y,t)=\gt e^{-\frac{\abs{y}^2}{4T}}\}}{}V(x,t)dS_\gt(x,t)\gt d\gt\\[4mm]
\phantom{\dint_{\!\!E\cap P_\gd}H(x-y,t)V(x,t)dxdt}
\leq\left[\gt\myint{0}{\gt}\myint{\{(x,t)\in Q_T^{B_R}:H(x-y,t)=\gs e^{-\frac{\abs{y}^2}{4T}}\}}{}V(x,t)dS_\gs(x,t) d\gs\right]_{\gt=0}^{\gt=\gd}
\\[4mm]
\phantom{\dint_{\!\!E\cap P_\gd}H(x-y,t)V(x,t)dxdt--}
-\myint{0}{\gd}\myint{0}{\gt}\myint{\{(x,t)\in Q_T^{B_R}:H(x-y,t)=\gs e^{-\frac{\abs{y}^2}{4T}}\}}{}V(x,t)dS_\gs(x,t) d\gs d\gt
\\[4mm]
\phantom{\dint_{\!\!E\cap P_\gd}H(x-y,t)V(x,t)dxdt}
\leq \gd\myint{0}{\gd}\myint{\{(x,t)\in Q_T^{B_R}:H(x-y,t)=\gs e^{-\frac{\abs{y}^2}{4T}}\}}{}V(x,t)dS_\gs(x,t) d\gs.
\EA$$
The first integration by parts is justified since $V\in L^1(Q_T^{B_R})$. Notice that
$$\BA {l}
 \gd\myint{0}{\gd}\myint{\{(x,t)\in Q_T^{B_R}:H(x-y,t)=\gs e^{-\frac{\abs{y}^2}{4T}}\}}{}V(x,t)dS_\gs(x,t) d\gs
 =\gd\dint_{\!\!Q_T^{B_R}\cap P_\gd}V(x,t)dxdt
\EA$$
and
$$\BA{ll}\gd\dint_{\!\!Q_T^{B_R}\cap P_\gd}V(x,t)dxdt\leq \gd\dint_{\!\!Q_T^{B_R}\cap R_\gd(y,T)}V(x,t)dxdt\\[4mm]
\phantom{\gd\dint_{\!\!Q_T^{B_R}\cap P_\gd}V(x,t)dxdt}
\leq\gb r^{-n}\myint{0}{\ga r}\myint{B_R\cap B_{(\ga r)^2}(y)}{}V(x,t)dxdt,
\EA$$
for some $\ga,\gb>0$ and if we have set $r=\gd^{-\frac{1}{n}}$. Notice also that $B_R\cap B_{(\ga r)^2}(y)=\emptyset$ if
$\abs{y}\geq R+(\ga r)^2$, or, equivalently, if $\abs{y}\geq R+\ga^2 \gd^{-\frac{2}{n}}$.\smallskip

(i) If $\abs{y}\geq R+\ga^2 $, we fix $\gd$ such that $1<\gd$, then
$$e^{\frac{\abs{y}^2}{4T}}\dint_{\!\!E}H(x-y,t)V(x,t)dxdt\leq \gd\dint_{\!\!E}V(x,t)dxdt,
$$
which can be made smaller than $\ge$ provided $\abs{E}$ is small enough.\smallskip

(ii) If $\abs{y}< R+\ga^2 $, then
$$\BA {l}e^{\frac{\abs{y}^2}{4T}}\dint_{\!\!E\cap P^c_\gd}H(x-y,t)V(x,t)dxdt
\leq e^{\frac{R^2+\ga^4}{2T}}\dint_{\!\!E\cap P^c_\gd}H(x-y,t)V(x,t)dxdt\\[4mm]
\phantom{e^{\frac{\abs{y}^2}{4T}}\dint_{\!\!E\cap P^c_\gd}H(x-y,t)V(x,t)dxdt}
\leq\gd  e^{\frac{R^2+\ga^4}{2T}}\dint_{\!\!E}V(x,t)dxdt.
\EA$$
Given $\ge>0$, we fix $\gd=r^{-n}$ such that
$$\BA {l}
e^{\frac{R^2+\ga^4}{2T}}\dint_{\!\!E\cap P_\gd}H(x-y,t)V(x,t)dxdt\leq
\gb e^{\frac{R^2+\ga^4}{2T}} r^{-n}\myint{0}{\ga r}\myint{B_R\cap B_{(\ga r)^2}(y)}{}V(x,t)dxdt\leq\myfrac{\ge}{2},
\EA$$
and then $\eta>0$ such that $\abs{E}\leq\eta$ implies
$$e^{\frac{\abs{y}^2}{4T}}\dint_{\!\!E\cap P^c_\gd}H(x-y,t)V(x,t)dxdt
\leq\gd  e^{\frac{R^2+\ga^4}{2T}}\dint_{\!\!E}V(x,t)dxdt\leq\myfrac{\ge}{2}.
$$
Therefore
$$e^{\frac{\abs{y}^2}{4T}}\dint_{\!\!E}H(x-y,t)V(x,t)dxdt\leq \ge,
$$
which is (\ref{stabV}).\hfill$\Box$\medskip

\noindent{\bf Remark} In Theorem \ref{uncondi} and Proposition \ref{suff}, the assumption of uniformity with respect to $y\in\BBR^n$ in (\ref{stabV}), (\ref{stabV+}) and  (\ref{stabV'}) can be replaced by uniformity with respect to $y\in B_{R_0}$ if all the measures $\gm_k$ have their support in $B_{R_0}$. A extension of these assumptions, valid when the convergent measures $\gm_k$ have their support in a fixed compact set is to assume that $V$ is {\bf locally strongly subcritical}, which means that (\ref{stabV}) holds uniformly with respect to $y$ in a compact set. Similar extension holds for (\ref{stabV'}).

\section{The supercritical case}
\setcounter{equation}{0}
\subsection{Capacities}
All the proofs in this subsection are similar to the ones of \cite{VY} and inspired by \cite{Fu};  we omit them. We assume also that there exists a positive measure $\xm_0$ such that $\mathbb{H}[\xm_0]V\in L^1(Q_T)$.

\begin{defin}
If $\xm\in\mathfrak{M}_+(\mathbb{R}^n)$ and $f$ is a nonnegative measurable function defined in $\xO$ such that
$$(t,x,y)\mapsto \mathbb{H}[\xm](y,t)V(x,t)f(x,t)\;\in \;L^1(Q_T\times\mathbb{R}^n;\;dxdt\otimes d\xm),$$
we set
$$\mathcal{E}(f,\xm)=\int_{Q_T}\left(\int_{\mathbb{R}^n} H(x-y,t)d\xm(y)\right)V(x,t)f(x,t)dxdt.$$
\end{defin}
If we put
$$\breve{\mathbb{H}}[f](y)=\int_{Q_T}H(x-y,t)V(x,t)f(x,t)dxdt,$$
then by Fubini's Theorem, $\breve{\mathbb{H}}[f](y)<\infty,$ $\xm-$almost everywhere in $\mathbb{R}^n$ and
$$\mathcal{E}(f,\xm)=\int_{\mathbb{R}^n}\left(\int_{Q_T} H(x-y,t)V(x,t)f(x,t)dxdt\right)d\xm(y).$$
\begin{prop}
Let $f$ be fixed. Then\\
(a) $y\mapsto\breve{\mathbb{H}}[f](y)$ is lower semicontinuous in $\mathbb{R}^n.$\\
(b) $\xm\mapsto\mathcal{E}(f,\xm)$ is lower semicontinuous in $\mathfrak M_+(\mathbb{R}^n)$ in the weak* topology.
\end{prop}

\begin{defin}
We denote by $\mathfrak{M}^V(\mathbb{R}^n)$ the set of all measures $\xm$ on $\mathbb{R}^n$ such that $V\mathbb{H}[|\xm|]\in L^1(Q_T).$ If $\xm$ is such a measure, we set
$$||\xm||_{\mathfrak {M}^V}=\int_{Q_T}\left(\int_{\mathbb{R}^n} H(x-y,t)d|\xm|(y)\right)V(x,t)dxdt=||V\mathbb{H}[|\xm|]||_{L^1(Q_T)}.$$
\end{defin}
If $E\subset\mathbb{R}^n$ is a Borel set, we put
$$\mathfrak {M}_+(E)=\{\xm\in\mathfrak {M}_+(\mathbb{R}^n):\;\xm(E^c)=0\}\quad\mathrm{and}\quad
\mathfrak{M}^V_+(E)=\mathfrak {M}^V(\mathbb{R}^n)\cap\mathfrak {M}_+(E).$$
\begin{defin}
If $E\subset\mathbb{R}^n$ is any borel subset we define the set function $C_V$ by
$$C_V(E):=\sup\{\xm(E):\;\xm\in\mathfrak M^V_+(E),\;||\xm||_{\mathfrak M^V}\leq1\};$$
this is equivalent to,
$$C_V(E):=\sup\left\{\frac{\xm(E)}{||\xm||_{\mathfrak M^V}}:\;\xm\in\mathfrak M^V_+(E)\right\}.$$
\end{defin}
\begin{prop}
The set function $C_V$ satisfies
$$C_V(E)\leq\sup_{y\in E}\left(\int_{Q_T}H(x-y,t)V(x,t)dxdt\right)^{-1}\quad\forall\;E\subset\mathbb{R}^n,\;E\;\mathrm{Borel}.$$
Furthermore equality holds if $E$ is compact. Finally,
$$C_V(E_1\cup E_2)=\sup\{C_V(E_1),C_V(E_2)\}\quad\forall\;E_i\subset\mathbb{R}^n,\;E_i\;\mathrm{Borel}.$$
\end{prop}
\begin{defin}
For any Borel $E\subset\mathbb{R}^n,$ we set
$$C_V^*(E):=\inf\{||f||_{L^\infty}:\;\breve{\mathbb{H}}[f](y)\geq1\;\forall y\in E\}.$$
\end{defin}
\begin{prop}
For any compact set $E\subset\mathbb{R}^n,$
$$C_V^*(E)=C_V(E).$$
\end{prop}

\subsection{The singular set of $V$}
In this section we assume that $V$ satisfies (\ref{equ16}), although much weaker assumption could have been possible. We define the singular set of $V$, $Z_V$ by
\begin{equation}\label{singV}
Z_V=\left\{x\in\mathbb{R}^n:\dint_{\!\!Q_T}H(x-y,t)V(y,t)dydt=\infty\right\}.
\end{equation}
Since the function $x\mapsto f(x)=\dint_{\!\!Q_T}H(x-y,t)V(y,t)dydt$ is lower semicontinuous, it is a Borel function and $Z_V$ is a Borel set.
\begin{lemma}
If $x\in Z_V$ then for any $r>0,$
$$\dint_{\!\!Q^{B_r(x)}_T}H(x-y,t)V(y,t)dydt=\infty.$$
\end{lemma}
\Proof We will prove it by contradiction, assuming that there exists  $r>0$, such that
$$\dint_{\!\!Q^{B_r(x)}_T}H(x-y,t)V(y,t)dy\leq M.$$
Replacing $H$ by its value, we derive
\bea
\nonumber
\dint_{\!\!Q_T}H(x-y,t)V(y,t)dydt&=&\dint_{\!\!Q^{B_r(x)}_T}H(x-y,t)V(y,t)dydt+\dint_{\!\!Q^{B^c_r(x)}_T}H(x-y,t)V(y,t)dydt\\ \nonumber
&\leq& M+C(n)\int_0^Tt^{-\frac{n+2}{2}}e^{-\frac{r^2}{4t}}dt<\infty.
\eea
Which is clearly a contradiction.\hfill$\Box$
\begin{lemma}
If $\xm$ is an admissible positive measure then $\xm(Z_V)=0.$\label{adm}
\end{lemma}
\Proof
Let $K\subset Z_V$ be a compact set. In view of the above lemma there exists a $R>0$ such that $K\subset B_R$ and
for each $x\in K ,$ we have
\be
\dint_{\!\!Q_T^{B_{2R}}}H(x-y,t)V(y)dy=\infty\label{par1}
\ee
and
\be
\dint_{\!\!Q_T^{B^c_{2R}}}H(x-y,t)V(y)dy<\infty.\label{par2}
\ee
Now, $\xm_K=\chi_{K}\xm$ is an admissible measure and by Fubini theorem we have
$$\BA{lll}
\displaystyle
\dint_{\!\!Q_T}\left(\int_{\mathbb{R}^n} H(x-y,t)d\xm_K(y)\right)V(x,t)dxdt=
\displaystyle\int_K \dint_{\!\!Q_T}H(x-y,t) V(x,t)dxdtd\xm(y)\\ [4mm]
\phantom{-\dint_{\!\!Q_T}\left(\int_{\mathbb{R}^n} H(x-y,t)d\xm_K(y)\right)V(x,t)dxdt}
=\displaystyle \int_K\dint_{\!\!Q_T^{B_{2R}}}H(x-y,t)V(x)dxdtd\xm(y)\\ [4mm]
\phantom{---\int_{Q_T}\left(\int_{\mathbb{R}^n} H(x-y,t)d\xm_K(y)\right)V(x,t)dxdt}+
\displaystyle\int_K\dint_{\!\!Q_T^{B^c_{2R}}}H(x-y,t)V(x)dxdtd\xm{y}.
\EA$$
By (\ref{par2}) the second integral above is finite and by (\ref{par1})
$$\dint_{\!\!Q_T^{B_{2R}}}H(x-y,t)V(x)dxdt=\infty\qquad\forall y\in K.$$
It follows that $\xm(K)=0.$ This implies $\xm(Z_V)=0$ by regularity.\hfill$\Box$
\begin{theorem}
If $\xm\in\mathfrak M_T(\mathbb{R}^n),$ $\gm\geq 0$ such that $\xm(Z_V)=0$, then $\xm$ is a good measure.
\end{theorem}
\Proof
We set $\xm_R=\chi_{B_R}\xm.$ By Proposition \ref{goodbound}, since $Z_V^{B_R}\subset Z_V$, $\xm_R$ is a good measure in $B_R$ with corresponding solution $u_\xm^R.$ In view of Lemma \ref{fragma}, $u_\xm^R$ satisfies
$$\displaystyle \dint_{\!\!Q_T^{B_{R}}}|u_\xm^R|\frac{n}{4(T-t)}e^{-\frac{|x|^2}{4(T-t)}}dxdt+\dint_{\!\!Q_T^{B_{R}}} |u_\xm^R|Ve^{-\frac{|x|^2}{4(T-t)}} dxdt\leq \int_{B_R}e^{-\frac{|x|^2}{4T}}d\xm.$$
Also $\{u_\xm^R\}$ is an increasing function, thus converges to $u_\xm.$ By the above estimate we have that $u_\gm$ belong to class $\mathcal{E_V}(Q_T)$ and is a weak solution of (\ref{sub1}).\hfill$\Box$
\begin{prop}
Let $\xm\in\mathfrak M_+(\mathbb{R}^n).$ Then $\xm(Z_V)=0$ if and only if there exists an increasing sequence of positive admissible measures which converges to $\xm$ in the weak* topology.
\end{prop}
\Proof The proof is similar as the one of \cite[Th 3.11]{VY} and we present it for the sake of completeness.
First, we assume that $\xm(Z_V)=0.$ Then we define the set
$$K_N=\left\{x\in\mathbb{R}^n:\int_{Q_T}H(x-y,t)V(y)dydt\leq N\right\}.$$
We note that $Z_V\cap K_N=\emptyset.$
We set $\xm_n=\chi_{K_N}\xm$ then we have
$$\int_{Q_T}\left(\int_{\mathbb{R}^n} H(x-y,t)d\xm_n(y)\right)V(x,t)dxdt\leq \xm(K_N).$$
Thus $\xm_n$ is admissible, increasing with respect $n.$ By the monotone theorem it follows that $\xm_n\to\chi_{Z_V^c}\xm.$ Since $\xm(Z_V)=0$ the result follows in this direction.\\
For the other direction. Let $\{\xm_n\}$ be an increasing sequence of positive admissible measure. Then by Lemma \ref{adm} we have that $\xm_n(Z_V)=0,\;\forall\;n\geq1.$ Since $\xm_n\leq\xm,$ there exist an increasing functions $h_n$ $\xm-$integrable such that $\xm_n=h_n\xm.$ Since $0=\xm_n(Z_V)\to\xm(Z_V)$ the result follows.\hfill$\Box$
\subsection{Properties of positive solutions and representation formula}
We first recall the construction of the kernel function for the operator $w\mapsto \prt_tw-\Gd w+Vw$ in $Q_T$, always assuming that $V$ satisfies (\ref{equ16}). For $\gd>0$ and $\gm\in \mathfrak M_T$, we denote by $w_\gd$ the solution of
\begin{equation}\label{E1}
\BA{lll}
\prt_tw-\xD w+V_\gd w=0,\qquad &\mathrm{in }\;\; Q_{T}\\
\phantom{\prt_tw-V_\gd\xD }
w(.,0)=\gm\qquad &\mathrm{in }\;\;\BBR^n.
\EA\end{equation}
where $V_\gd=V\chi_{Q_{\gd,T}}$ and $Q_{\gd,T}=(\gd,T)\times\BBR^n$. Then
\begin{equation}\label{approx0}
w_\gd(x,t)=\int_{\BBR^n}H_{V_\gd}(x,y,t)d\gm(y)
\end{equation}

\begin{lemma}\label{kern}
The mapping $\gd\mapsto H_{V_\gd}(x,y,t)$ is increasing and converges to  $H_V\in C(\BBR^n\times\BBR^n\times (0,T]))$ when $\gd\to 0$. Furthermore
there exists a function $H_V\in C(\BBR^n\times\BBR^n\times (0,T]))$ such that for any $\gm\in \mathfrak M_T(\BBR^n)$
\begin{equation}\label{E2}
\lim_{\gd\to 0}w_\gd(x,t)=w(x,t)=\int_{\BBR^n}H_{V}(x,y,t)d\gm(y).
\BA{lll}
\EA\end{equation}
\end{lemma}
\Proof Without loss of generality we can assume $\gm\geq 0$. By the maximum principle $\gd\mapsto H_{V_\gd}(x,y,t)$ is increasing and the result follows by the monotone convergence theorem.\hfill$\Box$\medskip

\noindent  If $\BBR^n$ is replaced by a smooth bounded domain $\Gw$, we can consider the problem
\begin{equation}\label{E1'}
\BA{lll}
\prt_tw-\xD w+V_\gd w=0\qquad &\mathrm{in }\;\; Q^\Gw_{T}\\
\phantom{\prt_t-\xD w+V_\gd w}
w=0\qquad &\mathrm{in }\;\;  \partial_lQ_T^\xO:=\prt\Gw\times (0,T] \\
\phantom{\prt_tw-V_\gd\xD w}
w(.,0)=\gm\qquad &\mathrm{in }\;\;\Gw.
\EA\end{equation}
where $V'_\gd=V\chi_{Q^\Gw_{\gd,T}}$ and $Q^\Gw_{\gd,T}=(\gd,T)\times\Gw$. Then
\begin{equation}\label{approx1}
w_\gd(x,t)=\int_{\Gw}H^\Gw_{V_\gd}(x,y,t)d\gm(y)
\end{equation}

The proof of the next result is straightforward.
\begin{lemma}\label{kern'}
The mapping $\gd\mapsto H^\Gw_{V_\gd}(x,y,t)$ increases and converges to  $H^\Gw_V\in C(\Gw\times\Gw\times (0,T]))$ when $\gd\to 0$. Furthermore
There exists a fonction $H^\Gw_V\in C(\Gw\times\Gw\times (0,T]))$ such that for any $\gm\in \mathfrak M_b(\Gw)$
\begin{equation}\label{E3}
\lim_{\gd\to 0}w_\gd(x,t)=w(x,t)=\int_{\Gw}H^\Gw_{V}(x,y,t)d\gm(y).
\BA{lll}
\EA\end{equation}
Furthermore $H^\Gw_V\leq H^{\Gw'}_V\leq H_V$ if $\Gw\subset\Gw'$.
\end{lemma}

\medskip

It is important to notice that the above results do not imply that $w$ is a weak solution of problem (\ref{equ1}). This question will be considered later on with the notion of reduced measure.
\begin{lemma}\label{lem1}
Assume $\xm\in\mathfrak M_+(\mathbb{R}^n)$ is a good measure and let $u$ be a positive weak solution of problem (\ref{sub1}). If $\xO$ is a smooth bounded domain, then there exists a unique positive weak solution $v$ of problem
\begin{equation}\label{prob*}
\BA{lll}
\prt_tv-\xD v+Vv=0,\qquad &\mathrm{in }\;\; Q_T^\xO,\\
\phantom{\prt_tv-\xD +Vv}
v=0\qquad&\mathrm{on }\;\; \partial_lQ_T^\xO\\
\phantom{\prt_tv-\xD v}
v(.,0)=\chi_{\xO}\xm\qquad &\mathrm{in }\;\;\xO.
\EA\end{equation}
Furthermore
\begin{equation}\label{E3'}
v(x,t)=\int_{\Gw}H^\Gw_{V}(x,y,t)d\gm(y).
\end{equation}
\end{lemma}
\Proof
Let $\{t_j\}_{j=1}^\infty$ be a sequence decreasing  to $0$, such that $t_j<T,\;\forall
\;j\in\mathbb{N}.$ We consider the following problem
\begin{equation}\label{S1}
\BA{lll}
\prt_tv-\xD v+Vv=0,\qquad &\mathrm{in}  &\xO\times(t_j,T],\\
\phantom{\prt_t-\xD v+Vv}
v=0\qquad&\mathrm{on} &\partial\xO\times(t_j,T]\\
\phantom{\prt_t-\xD v,}
v(.,t_j)=u(.,t_j)\qquad &\mathrm{in}&\xO\times\{t_j\}.
\EA\end{equation}
Since $u,Vu\in L^1(Q_T^{B_R})$ for any $R>0$, $t\mapsto u(.,t)$ is continuous with value in $L^1_{loc}(\BBR^n)$, therefore $u(.,t_j)\in L^1_{loc}(\BBR^n)$ and there exists a unique solution $v_j$ to (\ref{S1}) (notice also that $V\in L^\infty(Q_T^{B_R})$). By the maximum principle $0\leq v_j\leq u$ and by  standard parabolic estimates,  we may assume that the sequence $v_j$ converges locally uniformly in $\xO\times(0,T]$ to a function $v\leq u.$
Also, if $\xf\in C^{1,1;1}(\overline{Q_T^\xO})$ vanishes on $\partial_lQ_T^\xO$ and satisfies $\xf(x,T)=0,$ we have
\[
-\int_{t_j}^T\int_\xO v_j(\prt_t\xf+\xD\xf)dxdt+\int_{t_j}^T\int_\xO Vv_j\xf dxdt
+\int_\xO\xf(x,T-t_j)v_j(x,T-t_j)dx
=\int_\xO\xf(x,0)u(x,t_j)dx,
\]
where in the above equality we have taken $\xf(.,.-t_j)$ as test function. Since $\xf(.,T-t_j)\to 0$ uniformly and $u(.,t_j)\to\gm$ in the weak sense of measures, it follows by the dominated convergence theorem that
\[
-\dint_{\!\!Q_T^\Gw} v(\prt_t\xf+\xD\xf)dxdt+\dint_{\!\!Q_T^\Gw}Vv\xf dxdt
=\int_\xO\xf(y,0)d\gm(y),
\]
thus $v$ is a weak solution of problem (\ref{prob*}). Uniqueness follows as in Lemma \ref{monad}. Finally, for $\gd>0$, we consider the solution $w_\gd$ of (\ref{E1'}).  Then it is expressed by (\ref{approx0}). Furthermore
\[
-\dint_{\!\!Q^\Gw_T} w_\gd(\prt_t\xf+\xD\xf)dxdt+\dint_{\!\!Q^\Gw_T}  V_\gd w_\gd\xf dxdt
=\int_\xO\xf(x,0)d\gm(x),
\]
The sequence $w_\gd$ is decreasing, with limit $w$. Since $w_\gd\geq v$, then $w\geq v$. If we assume $\gf\geq 0$, it follows from dominated convergence and Fatou's lemma that
\[
-\dint_{\!\!Q^\Gw_T} w(\prt_t\xf+\xD\xf)dxdt+\dint_{\!\!Q^\Gw_T}  Vw\xf dxdt
\leq \int_\xO\xf(x,0)d\gm(x),
\]
Thus $w$ is a subsolution for problem (\ref{prob*}) for which we have comparison when existence. Finally $w=v$ and
(\ref{E3'}) holds.
\hfill$\Box$
\begin{lemma}
Assume $\xm\in\mathfrak M_+(\mathbb{R}^n)$ is a good measure and let $u$ be a positive weak solution of problem (\ref{sub1}).  Then for any  $(x,t)\in\mathbb{R}^n\times(0,T],$ we have
$$\lim_{R\to\infty}u_R= u,$$
where $\{u_R\}$ is the increasing sequence of the weak solutions of the problem (\ref{prob*}) with $\xO=B_R.$ Moreover, the convergence is uniform in any compact subset of $\mathbb{R}^n\times(0,T]$ and we have the representation formula
$$u(x,t)=\int_{\mathbb{R}^n}H_V(x,y,t)d\xm(y).$$\label{repr1}
\end{lemma}
\Proof By the maximum principle,  $u_R\leq u_{R'}\leq u$ for any $0<R\leq R'.$ Thus $u_R\to w\leq u.$ Also by standard parabolic estimates, this convergence is locally uniformly. Now by dominated convergence theorem, it follows that $w$ is a weak solution of problem (\ref{sub1}) with initial data $\xm.$
Now we set $\widetilde{w}=u-w\geq0.$ Since $\widetilde{w}$ satisfies in the weak sense

$$\BA {lll}
\widetilde{w}_t-\xD \widetilde{w}+V\widetilde{w}= 0\qquad &\mathrm{in }\;  Q_T\\[1mm]
\phantom{\widetilde{w}_t+V\widetilde{w}}
\widetilde{w}(x,t)\geq 0 \qquad &\mathrm{in }\; Q_T
\\[1mm] \phantom{\widetilde{w}_t+V\widetilde{w}}
\widetilde{w}(x,0)=0\qquad &\mathrm{in }\; \mathbb{R}^n,
\EA$$
and $V\geq 0$, it clearly satisfies
$$\BA {lll}
\widetilde{w}_t-\xD \widetilde{w}\leq 0 \qquad &\mathrm{in }\;  Q_T\\[1mm]
\phantom{\widetilde{w}_t}
\widetilde w(x,t)\geq 0 \qquad &\mathrm{in }\; Q_T
\\[1mm] \phantom{\widetilde{w}_t}
\widetilde{w}(x,0)=0\qquad &\mathrm{in }\; \mathbb{R}^n,
\EA$$
which implies $\widetilde{w}=0.$ By the previous lemme $u_R$ admits the representation
\[
u^R(x,t)=\int_{B_R}H_V^{B_R}(x,y,t)d\xm(y).
\]
Since $\{H_V^{B_R}\}$ is an increasing sequence and $\lim_{R\to\infty}H_V^{B_R}=H_V$,  we have using again Fatou's lemma as in the proof of Lemma \ref{repr1}
$$u(x,t)=\lim_{R\to\infty}u^R(x,t)=\lim_{R\to\infty}\int_{B_R}H_V^{B_R}(x,y,t)d\xm(y)=\int_{\mathbb{R}^n}H_V(x,y,t)d\xm(y)$$
\hfill$\Box$
\begin{lemma}{\textbf{Harnack inequality}}
Let $C_1>0$ and $V(x,t)$ be a potential satisfying (\ref{equ16})
If $u$ is a positive solution of (\ref{equ17}),
then the Harnack inequality is valid:
$$u(y,s)\leq u(x,t)\exp \left(C(n,C_1)\left({\frac{|x-y|^2}{t-s}+\frac{t}{s}}+1\right)\right),\qquad\forall\;(y,s),(x,t)
\in Q_T\,,\;s<t.$$
\end{lemma}
\Proof
We extend $V$ for $t\geq T$ by the value $C_1t^{-1}$. We consider the linear parabolic problem
\be
\prt_tu-\xD u+Vu=0,\qquad \mathrm{in}\qquad \mathbb{R}^n\times[1,\infty),\label{2}
\ee
It is well known that, under the assumption (\ref{equ16}), every positive solution $u(x,t)$ of (\ref{2}) satisfies the Harnack inequality
$$
u(y,s)\leq u(x,t)\exp \left(C(n,C_1)\left({\frac{|x-y|^2}{t-s}+\frac{t}{s}}+1\right)\right),\qquad\forall\;(x,t)\in\mathbb{R}^n\times[1,\infty).
$$
Set $\tilde{u}(x,t)=u(\frac{t}{\xl^2}\frac{x}{\xl}).$ Then $\tilde{u}$ satisfies
$$\mathrm{u}_t-\xD \mathrm{u}+\frac{1}{\xl^2}V(\frac{t}{\xl^2}\frac{x}{\xl})\tilde{u}=0,\qquad \mathrm{in} \; \mathbb{R}^n\times(0,\infty).$$
We note here that $\frac{1}{\xl^2}V(\frac{t}{\xl^2}\frac{x}{\xl})\leq C_1,$ $\forall t\geq\frac{1}{\xl^2},$ thus
$\tilde{u}$ satisfies the Harnack inequality
$$\widetilde{u}(y,s)\leq \widetilde{u}(x,t)\exp \left(C(n,C_1)\left({\frac{|x-y|^2}{t-s}+\frac{t}{s}}+1\right)\right),\qquad\forall\;(x,t)
\in\mathbb{R}^n\times[\frac{1}{\xl^2},\infty).$$
By the last inequality and the definition of $\tilde{u}$ we derive the desired result.\hfill$\Box$\medskip

Next, we set
\begin{equation}\label{Sing}
\mathcal{S}ing_V(\mathbb{R}^n):=\{y\in\mathbb{R}^n:\;H_V(x,y,t)=0\}
\end{equation}
If $H_V(x,y,t)=0$ for some $(x,t)\in Q_T$, then $H_V(x',y,t')=0$ for any $(x',t')\in Q_T$, $t'<t$ by Harnack inequality principle. We prove the {\textbf{Representation formula}}.
\begin{theorem}\label{present}
Let $u$ be a positive solution of (\ref{equ17}).
Then there exists a measure $\xm\in\mathfrak M_+(\mathbb{R}^n)$ such that
$$u(x,t)=\int_{\mathbb{R}^n}H_V(x,y,t)d\xm(y),$$
and $\gm$ is concentrated on $(\mathcal{S}ing_V(\mathbb{R}^n))^c$.
\end{theorem}
\Proof By Lemma \ref{repr1}
we have
$$u(x,t)=\int_{\mathbb{R}^n}H_V(x,y,t-s)u(y,s)dy\;\;\;\mathrm{for\; any}\; s<t\leq T.$$
We assume that $s\leq\frac{T}{2}.$ By Harnack inequality on $x\mapsto H_V(x,y,\frac{T}{2})$
$$\int_{\mathbb{R}^n}H_V(0,y,\frac{T}{2})u(y,s)dy\leq c(n)\int_{\mathbb{R}^n}H_V(0,y,T-s,)u(y,s)dy=c(n)u(0,T).$$
For any Borel set $E,$ we define the measure $\gr_s$ by
$$\xr_s(E):
=\int_{E}H_V(0,y,\frac{T}{2})u(y,s)dy\leq\int_{\mathbb{R}^n}H_V(0,y,\frac{T}{2})u(y,s)dy\leq c(n)u(T,0).$$
Thus there exists a decreasing sequence $\{s_j\}_{j=1}^\infty$ which converges to origin, such that
the measure $\xr_{s_j}$ converges in the weak* topology to a positive Radon measure $\xr.$
Also we have the estimate $\xr(\mathbb{R}^n)\leq C(n)u(0,T).$
Now choose $(x,t)\in Q_T$ and $j_0$ large enough such that $t> s_{j_0}.$ Let $\xe>0,$ we set for any $j\geq j_0$,
$$W_j(y)=\frac{H_V(x,y,t-s_j)}{H_V(0,y,\frac{T}{2})+\xe}.$$
For any $R>0$ and $|y|>R$ we have
$$W_j(y)\leq\frac{1}{\xe}H_V(x,y,t-s_j)\leq \frac{1}{\xe}H(x-y,t-s_j)<\frac{1}{\xe}C(x,R,t-s_j),$$
where $\lim_{R\to\infty}C(x,R,t-s_j)=0.$ We have also
$$\int_{|y|\geq R}W_j(y)d\xr_j\leq \frac{1}{\xe}C(x,R,t-s_j)c(n)u(T,0).$$
For any $|y|<R$, we have by standard parabolic estimates that $W_j(y)\to \frac{H_V(x,y,t)}{H_V(\frac{T}{2},0,y)+\xe}$ when $j\to\infty$, uniformly with respect to $y.$
Thus by the above estimates it follows
$$\int_{\mathbb{R}^n}W_j(y)d\xr_j\to\int_{\mathbb{R}^n}\frac{H_V(x,y,t)}{H_V(0,y,\frac{T}{2})+\xe}d\xr.$$
For sufficiently large $j$ we have
$$\BA {ll}
\displaystyle
\int_{\mathbb{R}^n}\frac{H_V(x,y,t-s_j)}{H_V(0,y,\frac{T}{2})+\xe}d\xr_{s_j}=
\int_{\mathbb{R}^n}\frac{H_V(x,y,t-s_j)}{H_V(0,y,\frac{T}{2})+\xe}\left( H_V(0,y,\frac{T}{2})+\xe-\xe\right)u(y,s_j)dy\\
\phantom{\displaystyle\int_{\mathbb{R}^n}\frac{H_V(x,y,t-s_j)}{H_V(0,y,\frac{T}{2})+\xe}d\xr_{s_j}}
\displaystyle
=u(x,t)-\xe\int_{\mathbb{R}^n}\frac{ H_V(x,y,t-s_j,)}{H_V(0,y,\frac{T}{2})+\xe}u(y,s_j)dy.
\EA$$
Note that this is a consequence of the identity
$$\int_{\mathbb{R}^n}H_V(x,y,t-s_j)u(y,s_j)dy=u(x,t).$$
Thus as before, we define
$d\widetilde{\xr}_j=H_V(x,y,t-s_j)u(y,s_j)dy$ and thus there exists a subsequence, say $\{\widetilde{\xr}_j\}$ converges in the weak* topology to a positive Radon measure $\widetilde{\xr}.$ Thus we have
$$\BA {ll}
\displaystyle
\xe\int_{\mathbb{R}^n}\frac{ H_V(x,y,t-s_j)}{H_V(0,y,\frac{T}{2})+\xe}u(y,s_j)dy=
\xe\int_{\mathbb{R}^n}\chi_{\left(\mathcal{S}ing_V(\mathbb{R}^n)\right)^c}\frac{ H_V(x,y,t-s_j)}{H_V(0,y,\frac{T}{2})+\xe}u(s_j,y)dy\\
\phantom{\displaystyle \xe\int_{\mathbb{R}^n}\frac{ H_V(x,y,t-s_j)}{H_V(0,y,\frac{T}{2})+\xe}u(y,s_j)dy}
\displaystyle \to \xe\int_{\mathbb{R}^n}\chi_{\left(\mathcal{S}ing_V(\mathbb{R}^n)\right)^c}
\frac{1}{H_V(\frac{T}{2},0,y)+\xe}d\widetilde{\xr}.
\EA$$
Combining the above relations, we derive
\be
\int_{\mathbb{R}^n}\frac{H_V(x,y,t)}{H_V(0,y,\frac{T}{2})+\xe}d\xr=u(x,t)-
\xe\int_{\mathbb{R}^n}\chi_{\left(\mathcal{S}ing_V(\mathbb{R}^n)\right)^c}
\frac{1}{H_V(0,y,\frac{T}{2})+\xe}d\widetilde{\xr}.\label{3}
\ee
Now, we have
$$\lim_{\xe\to0}\chi_{\left(\mathcal{S}ing_V(\mathbb{R}^n)\right)^c}\frac{\xe}{H_V(\frac{T}{2},0,y)+\xe}=0,$$
and by Harnack inequality on the function $x\mapsto H_V(x,y,t)$
$$\frac{H_V(x,y,t)}{H_V(0,y,\frac{T}{2})+\xe}\leq C(t,T),$$
thus by dominated convergence theorem, we can let $\xe$ tend to $0$ in (\ref{3}) and obtain
$$\int_{\mathbb{R}^n}\frac{H_V(x,y,t)}{H_V(0,y,\frac{T}{2})}d\xr=u(x,t).$$
And the result follows if we set
$$d\xm=\chi_{\left(\mathcal{S}ing_V(\mathbb{R}^n)\right)^c}\frac{1}{H_V(0,y,\frac{T}{2})}d\xr.$$\hfill$\Box$

In the next result we give a construction of $H_V$, with some estimates and a different proof of the existence of an initial measure for positive solutions of (\ref{equ16}).

\begin{theorem} Assume $V$ satisfies (\ref{equ16}) and $u$ is a positive solution of (\ref{equ17}) then there exists a positive Radon measure $\xm$ in $\mathbb{R}^n$ such that
\be\label{R1}
u(x,t)=\int_{\mathbb{R}^n}\e^{\psi (x,t)}\Gamma (x,y,t,0)d\xm(y)
\ee
where
\be\label{R2}
\psi (x,t)=\int_t^T\int_{\BBR^n}\myfrac{e^{-\frac{\abs{x-y}^2}{4(s-t)}}}{4\gp(t-s)}V(y,s)dyds
\ee
and
\be\label{R3}
c_1\myfrac{e^{-\gamma_1\frac{\abs{x-y}^2}{s-t}}}{(t-s)^\frac{n}{2}}\leq \Gamma (x,y,t,s)\leq c_2\myfrac{e^{-\gamma_2\frac{\abs{x-y}^2}{s-t}}}{(t-s)^\frac{n}{2}}
\ee
for some positive constants $c_i$ and $\gamma_i$, $i=1,2$.\label{111*}
\end{theorem}
\Proof Assuming that $u$ is a positive solution of (\ref{equ17}), we set $u(x,t)=e^{\psi(x,t)}v(x,t)$. Then
\be\label{R4}
\prt_t v-\Gd v-2\nabla\psi.\nabla v-\abs{\nabla\psi}^2v+(\prt_t\psi-\Gd\psi+V)v=0.
\ee
We choose $\psi$ as the solution of problem
\be\label{R5}\BA {ll}
-\prt_t\psi-\Gd\psi+V\psi=0\quad&\text{in }\; Q_T\\[2mm]
\phantom{-\prt_t\psi-\Gd\psi}
\psi(.,T)=0\quad&\text{in }\; \BBR^n.
\EA\ee
Then $\psi$ is expressed by (\ref{R2}). Furthermore, by standard computations,
\be\label{R6}\BA {ll}
(i)&\qquad 0\leq \psi(x,t)\leq c\ln\frac{T}{t}\phantom{-\prt_t\psi-\Gd\psi}\phantom{-\prt_t\psi-\Gd\psi}\\[2mm]
(ii)&\qquad \abs{\nabla \psi(x,t)}\leq c_1(T)+c_2(T)\ln\frac{T}{t}\phantom{-\prt_t\psi-\Gd\psi}\phantom{-\prt_t\psi-\Gd\psi}
\EA\ee
The function $v$ satisfies
\be\label{R7}
\prt_t v-\Gd v-2\nabla\psi.\nabla v-\abs{\nabla\psi}^2v=0.
\ee
Then, by (\ref{R6}),
\be\label{R8}\BA {ll}
(i)&\qquad 0\leq \displaystyle\int_{\BBR^n}\sup\{\abs{\psi(x,s)}^q:x\in\BBR^n\}ds\leq M_1\\[4mm]
(ii)&\qquad 0\leq \displaystyle\int_{\BBR^n}\sup\{\abs{\nabla \psi(x,s)}^q:x\in\BBR^n\}ds\leq M_2
\EA\ee
for any $1\leq q<\infty$ for some $M_i\in\BBR_+$. This is the condition $H$ in \cite{Ar} with $R_0=\infty$ and $p=\infty$. Therefore there exists a kernel function $\Gamma\in C(\BBR^n\times\BBR^n\times(0,T)\times (0,T))$ which satisfies
(\ref{R3}) and there exists also a positive Radon measure $\gm$ in $\BBR^n$ such that
\be\label{R9}\BA {ll}
v(x,t)=\displaystyle\int_{\BBR^n}\Gamma (x,y,t,0) d\gm (y).
\EA\ee
Finally $u$ verifies
 \be\label{R10}\BA {ll}
u(x,t)=e^{\psi(x,t)}\displaystyle\int_{\BBR^n}\Gamma (x,y,t,0) d\gm (y).
\EA\ee
\hfill$\Box$\\
We recall that $\mathcal{S}ing_V(\mathbb{R}^n):=\{y\in\mathbb{R}^n:\;H_V(x,y,t)=0\}.$
\begin{theorem}
Let  $\xd_\xi$ be the Dirac measure concentrated at $y$ and let $V$ satisfies (\ref{equ16}). Then
$$H_V(x,\xi,t)=\int_{\mathbb{R}^n}e^{\psi(x,t)}\xG(x,y,t)d\xm_\xi(y),$$
where $\xm_\xi$ is a positive Radon measure such that $$\xd_\xi\geq\xm_\xi,$$
and $\psi, \;\xG$ are the functions in (\ref{R2}) and (\ref{R3}) respectively.\\
Furthermore, if
$$\limsup_{t\rightarrow0}\psi(\xi,t)=\limsup_{t\rightarrow0}\myint{t}{T}\myint{\BBR^n}{}\left(\myfrac{1}{4\gp (s-t)}\right)^{\frac{n}{2}}e^{-\frac{\abs{\xi-y}^2}{4 (s-t)}}V(y,s) dy ds=\infty$$
then $$\xi\in\CS ing_V,\; \mathrm{i.e.} \;H_V(x,\xi,t)=0,\;\forall (x,t)\in\mathbb{R}^n\times(0,\infty). $$\label{thF}
\end{theorem}
{\it proof.} First we note that $H_{V_k}(x,\xi,t)$ is the solution of problem (\ref{equ20_*}) with $\xd_\xi$ as initial data. Since $H_{V_k}(x,\xi,t)\downarrow H_{V}(x,\xi,t),$ we have by maximum principle , $H(x,\xi,t)\geq H_V(x,\xi,t).$
Now by Theorem \ref{111*}, there exists a positive Radon measure $\xm_\xi$ in $\mathbb{R}^n$ such that
\be\label{111}
H_V(x,,\xi,t)=\int_{\mathbb{R}^n}\e^{\psi (x,t)}\Gamma (x,y,t,0)d\xm_\xi(y)
\ee
Let $\xf\in C_0(\mathbb{R}^n)$ then we have by the properties of $\xG(x,\xi,t)$ (see \cite{Ar}) and (\ref{111})
$$\lim_{t\rightarrow0}\int_{\mathbb{R}^n}H_V(x,\xi,t)\xf(x)dx\geq\lim_{t\rightarrow0}\int_{\mathbb{R}^n}\int_{\mathbb{R}^n}
\xG(x,y,t)\xf(x)dxd\xm_\xi(y)=\int_{\mathbb{R}^n}\xf(y)\xm_\xi(y),$$
That is
\be
\int_{\mathbb{R}^n}\xf(x)d\xd_\xi(x)\geq\int_{\mathbb{R}^n}\xf(x)d\xm_\xi(x)\Rightarrow\xd_\xi\geq\xm_\xi,\label{333}
\ee
since $\xf$ is an abstract function in space $C_0(\mathbb{R}^n).$\\
Also we have that there exist positive constants $C_1,\;C_2$ such that
\be
\xG(x,y,t)\geq C_1H(x,y,C_2 t). \label{222}
\ee
 Also we have
\bea
\nonumber
H(\xi,\xi,t)\geq H_V(\xi,\xi,t)&=&\int_{\mathbb{R}^n}H_V(\xi,y,t)d\xm_\xi(y)=\int_{\mathbb{R}^n}e^{\psi(\xi,t)}\xG(\xi,y,t)d\xm_\xi(y)\\ \nonumber
(\mathrm{by}\;\;(\ref{222}))\qquad&\geq& C_1\int_{B(\xi,\sqrt{C_2t})}e^{\psi(t,\xi)}H(\xi,y,C_2 t)d\xm_\xi(y)\\ \nonumber
(\mathrm{By\;Harnack \;inequality})\qquad&\geq& C(T,n,C_1,C_2)\int_{B(\xi,\sqrt{C_2t})}e^{\psi(\xi,t)}H(\xi,\xi,\frac{C_2t}{2})d\xm_\xi(y)\\ \nonumber
&=& C(T,n,C_1,C_2)e^{\psi(\xi,t)}H(\xi,\xi,\frac{C_2t}{2})\xm_\xi(B(\xi,\sqrt{C_2t}))
\eea
Thus by  the last inequality and the fact that $$\frac{H(\xi,\xi,t)}{H(\xi,\xi,\frac{C_2t}{2})}= C(C_2,n)>0,$$ we have
$$ C(T,n,C_1,C_2)\geq e^{\psi(t,\xi)}\xm_\xi(B(\xi,\sqrt{C_2t})).$$
But $limsup_{t\rightarrow0}\psi(\xi,t)=\infty$ which implies
$$\lim_{t\rightarrow0}\xm_\xi(B(\xi,\sqrt{C_2t})=\xm_\xi(\{\xi\})=0.$$
Thus by (\ref{333}) we have $\xm_\xi\equiv0,$ i.e. $H_V(x,\xi,t)=0,\;\forall (x,t)\in\mathbb{R}^n\times(0,\infty).$\hfill$\Box$
\subsection{Reduced measures}
In this section we assume that $V$ is nonnegative, but not necessarily satisfies (\ref{equ16}), therefore we can construct
$H_V[\gm]$ for
 $\xm\in\mathfrak M_{T}(\mathbb{R}^n)$. Furthermore, if $\gm$ is nonnegative we can consider the solution $u_k$
 of the problem
\begin{equation}\label{Red1}\BA {ll}
\prt_tu-\xD u+V^ku=0,\quad &\mathrm{in}\;\; Q_T\\ \phantom{\prt_tu-\xD u\!+}
u(.,0)=\xm\quad &\mathrm{in}\;\;\mathbb{R}^n,
\EA\end{equation}
where $V^k=\min\{V,k\}.$ Then there holds
$$u_k(x,t)=\int_{\mathbb{R}^n}H_{V^k}(t,x,y)d\xm(y)=\mathbb{H}_{V^k}[\xm](x,t),$$
and
$$u_k+\int_0^t\int_{\mathbb{R}^n}H(t-s,x,y)V^k u_kdyds=\mathbb{H}[\xm].$$
Since $k\mapsto H_{V^k}$ is decreasing and converges to $H_V,$ we derive
$$\lim_{k\to\infty}u_k=u=\int_{\mathbb{R}^n}H_{V}(t,x,y)d\xm(y).$$
By Fatou's lemma
$$\int_0^t\int_{\mathbb{R}^n}H(t-s,x,y)V udyds\leq\liminf_{k\to\infty}\int_0^t\int_{\mathbb{R}^n}H(t-s,x,y)V^k u_kdyds.$$
It follows
$$u(x,t)+\int_0^t\int_{\mathbb{R}^n}H(t-s,x,y)V udyds\leq\int_{\mathbb{R}^n}H_{V}(t,x,y)d\xm(y),\quad\forall (x,t)\in Q_T.$$
Now since $Vu\in L^1_{loc}(\overline{Q}_T)$ and
$$\prt_tu-\xD u+Vu=0,\qquad \mathrm{in}\qquad Q_T,$$
the function
$$u(x,t)+\int_0^t\int_{\mathbb{R}^n}H(t-s,x,y)V udyds$$
is nonnegative and satisfies the heat equation in $Q_T$. Therefore it admits an initial trace $\xm^*\in\mathfrak M_+(\mathbb{R}^n)$ and actually $\gm^*\in \mathfrak M_T(\BBR^n)$.
Furthermore, we have
$$u(x,t)+\int_0^t\int_{\mathbb{R}^n}H(t-s,x,y)V udyds=\int_{\mathbb{R}^n}H(x-y,t)d\xm^*(y),\;\forall (x,t)\in Q_T.,$$
or equivalently, $u$ is a positive weak solution of the problem
$$\BA {ll}
\nonumber
\prt_tu-\xD u+Vu=0\qquad &\mathrm{in }\;\; Q_T\\ \nonumber
\phantom{-\xD u--}
u(.,0)=\xm^* &\mathrm{in }\;\;\mathbb{R}^n.
\EA$$
Note that $\xm^*\leq\xm$ and the mapping $\xm\mapsto\xm^*$ is nondecreasing.
\begin{defin}
The measure $\xm^*$ is the reduced measure associated to $\xm$
\end{defin}
The proofs of the next two Propositions are similar to the ones of  \cite[Section 5]{VY}.
\begin{prop}
There holds $\mathbb{H}_V[\xm]=\mathbb{H}_V[\xm^*].$ Furthermore the reduced measure $\xm^*$ is the largest measure for which the following problem
\begin{equation}\BA {ll}\label{red1}
\prt_tv-\xD v+Vv=0\qquad &\mathrm{in }\;\; Q_T\\ [2mm]
\!\!\xl\in\mathfrak M_+(\mathbb{R}^n),\;\xl\leq\xm
\\ \phantom{----}
v(.,0)=\xl\qquad &\mathrm{in }\;\mathbb{R}^n,
\EA\end{equation}
admits a solution.
\end{prop}
\begin{prop}
Let $W_k$ be an increasing sequence of nonnegative bounded measurable functions converging to $V$ a.e. in $Q_T.$ Then the solution $u_k$ of
\begin{equation}\BA {ll}\label{red2}
\nonumber
\prt_tv-\xD v+W_kv=0\qquad &\mathrm{in }\;\; Q_T\\ \phantom{---,--}
v(.,0)=\xm\qquad &\mathrm{in}\;\;\mathbb{R}^n,
\EA\end{equation}
converges to $u_{\xm^*}.$
\end{prop}
We recall that $\mathcal{S}ing_V(\mathbb{R}^n):=\{y\in\mathbb{R}^n:\;H_V(x,y,t)=0\}.$
\begin{prop}
Let $\xm$ be a nonnegative measure in $\mathcal{M_T}(\mathbb{R}^n)$. Then\smallskip

\noindent (i) $(\gm-\gm^*)\left(\left(\mathcal{S}ing_V(\mathbb{R}^n)\right)^c\right)=0$

\noindent (ii) If $\xm\left(\left(\mathcal{S}ing_V(\mathbb{R}^n)\right)^c\right)=0,$ then $\xm^*=0.$\smallskip

\noindent (iii) There always holds $\mathcal{S}ing_V(\mathbb{R}^n)= Z_V.$
\end{prop}
{\it proof.} The proofs of (i), (ii) and the fact that $\mathcal{S}ing_V(\mathbb{R}^n)\subset Z_V$ are similar as in \cite[Section 5]{VY}, and we omit them. 

The proof of $Z_V\subset\mathcal{S}ing_V(\mathbb{R}^n)$ is a immediately consequence of Theorem \ref{thF}. Indeed, if $\xi\in Z_V$ then 
$$\limsup_{t\rightarrow0}\myint{t}{T}\myint{\BBR^n}{}\left(\myfrac{1}{4\gp (s-t)}\right)^{\frac{n}{2}}e^{-\frac{\abs{\xi-y}^2}{4 (s-t)}}V(y,s) dy ds=\infty,$$
thus $\xi\in\mathcal{S}ing_V(\mathbb{R}^n).$\hfill$\Box$
\section{Initial trace}
\subsection{The direct method}
The initial trace that we developed in this section is an adaptation to the parabolic case of the notion of boundary trace for elliptic equations (see \cite{MV},  \cite{MV1}, \cite{VY}). If $G\subset \overline{Q_T}$ is a relatively open set, we denote
$$W(G)=\displaystyle \bigcap_{1\leq p<\infty}W^{2,1}_p(G)\quad\text{and }\;W_{loc}(G)
=\displaystyle \bigcap_{1\leq p<\infty}W^{2,1}_{p\, loc}(G).$$
 Since $V\in L^{\infty}_{loc}(Q_T)$, any solution of (\ref{equ17}) belongs to $W_{loc}(Q_T)$.
\begin{prop}\label{s(u)}
Let $u\in W_{loc}(Q_T)$ be a positive solution (\ref{equ17}).
Assume that, for some $x\in\BBR^n,$ there exists an open bounded neighborhood $U$ of $x$ such that
\begin{equation}\label{tr1}
\dint_{\!\!Q^U_T}u(y,t)V(y,t)dxdt<\infty
\end{equation}
Then $u\in L^1(U\times(0,T))$ and there exists a unique positive Radon measure $\xm$ in $U$ such that
$$\lim_{t\to0}\int_{U}u(y,t)\xf(x)dx=\int_{U}\xf(x)d\xm,\;\;\;\;\forall\xf\in C_0^\infty(U).$$
\end{prop}
\Proof
Since $Vu\in L^1(U\times(0,T))$ the following problem has a weak solution $v$ (see \cite{MV}).
\bea
\nonumber
\prt_tv-\xD v&=&Vu,\qquad \mathrm{in}\qquad U\times(0,T],\\ \nonumber
v(x,t)&=&0\qquad\mathrm{on}\;\; \partial U\times(0,T]\\ \nonumber
v(x,0)&=&0\qquad \mathrm{in}\;\;U.
\eea
Thus the function $w=u+v$ satisfies the heat equation. Thus there exists a unique Radon measure $\xm$ such that
$$\lim_{t\to0}\int_{U}w(y,t)\xf(x)dx=\int_{U}\xf(x)d\xm,\;\;\;\;\forall\xf\in C_0^\infty(U).$$
And the result follows since the initial data of $v$ is zero.\hfill$\Box$\\

We set
\begin{equation}\label{tr2}
\mathcal{R}(u)=\left\{y\in\mathbb{R}^n:\;\exists\; \mathrm{bounded \;neighborhood}\; U\; \mathrm{of}\;y,\;\dint_{\!\!Q^U_T}u(y,t)V(y,t)dxdt<\infty\right\}.
\end{equation}
Then $\mathcal{R}(u)$ is open and there exists a unique positive Radon measure $\xm$ on $\mathcal{R}(u)$ such that
\be
\lim_{t\to 0}\int_{\mathcal{R}}u(y,t)\xf(x)dx=\int_{\mathcal{R}}\xf(x)d\xm,\;\;\;\;\forall\xf\in C_0^\infty(\mathcal{R}).\label{tracemeasure}
\ee
\begin{prop}
Let $u\in W_{loc}(\mathbb{R}^n\times(0,T])$ be a positive solution of (\ref{equ17}).
Assume that, for some $x\in\mathbb{R}^n,$ there holds
\begin{equation}\label{tr3}\dint_{\!\!Q^U_T}u(y,t)V(y,t)dydt=\infty
\end{equation}
 for any bounded open neighborhood $U$ of $x$. Then
\begin{equation}\label{singular}
\limsup_{t\to 0}\int_{U}u(y,t)dy=\infty.
\end{equation}
\end{prop}
\Proof
We will prove it by contradiction. We assume that there exists an open neighborhood of $x$ such that
$$\int_{U}u(y,t)dy\leq M<\infty\qquad\forall t\in(0,T).$$
Then $\norm u_{L^1(Q^U_T)}\leq MT$.
Let $B_r(x)\subset\subset U$ for some $r>0$, and $\xz\in C_0^\infty(B_r(x)),$ such that
$\xz=1$ in $B_{\frac{r}{2}}(x),$ $\xz=0$ in $B_r^c(x)$ and $0\leq\xz\leq1.$ Then since $u$ is a positive solution we have
$$\int_U\prt_tu\xz dx-\int_Uu\xD\xz dx+\int_UVu\xz dx=0\Rightarrow
\int_{B_{\frac{r}{2}}}Vu dx\leq\int_U\prt_tu\xz dx-\int_Uu\xD\xz dx\Rightarrow$$
$$\int_U\prt_tu\xz dx-\int_Uu\xD\xz dx+\int_UVu\xz
dx=0\Rightarrow
\int_{B_{\frac{r}{2}}}Vu dx\leq -\int_U\prt_tudx+ M\norm{\Gd\gz}_{L^\infty}.$$
Integrating the above inequality on $(s,T)$, we get
\begin{equation}\label{tr4}
\int_s^T\int_{B_{\frac{r}{2}}}Vu dxdr\leq -\int_Uu(x,T)dx+\int_Uu(x,s)dx+\norm u_{L^1(Q^U_T)}\norm{\Gd\gz}_{L^\infty}.
\end{equation}
Letting $s\to 0$, we reach a contradiction.\hfill$\Box$\medskip

\noindent {\bf Remark.} It is not clear wether there holds
\begin{equation}\label{singular1}
\liminf_{t\to 0}\int_{U}u(y,t)dy=\infty.
\end{equation}
However, it follows from (\ref{tr4}) that if $u\in L^1(Q^U_T)$, the above equality holds.\medskip

\begin{defin}
If $u$ is a positive solution of (\ref{equ17}), we set $\CS(u)=\BBR^n\setminus\CR(u)$. The couple $(\CS(u),\gm)$ is called the initial trace of $u$, denoted by $tr_{\{t=0\}}(u)$. The sets $\mathcal{R}(u)$ and $\mathcal{S}(u)$ are respectively the regular and the singular sets of $tr_{\{t=0\}}(u)$ and $\gm\in\mathfrak M_+(\CR(u))$ is its regular part.
\end{defin}

\noindent{\bf Example} Take $V(x,t)=ct^{-1}$, $c>0$. If $u$ satisfies
\begin{equation}\label{singular2}\prt_tu-\Gd u+\frac{c}{t}u=0
\end{equation}
then $v(x,t)=t^cu(x,t)$ satisfies the heat equation. Thus, if $u\geq 0$, there exists $\gm\in\mathfrak M_+(\BBR^n)$ such that
\begin{equation}\label{singular3}
u(x,t)=t^{-c}\BBH[\gm](x,t).
\end{equation}
This is a representation formula. Notice that $Vu(x,t)=ct^{-c-1}\BBH[\gm](x,t)$, therefore the regular set of $tr_{\{t=0\}}(u)$ may be empty.

\begin{prop}\label{subsol}
Assume $V$ satisfies (\ref{equ16}) and let $u\in W_{loc}(Q_T)$ be a positive solution of (\ref{equ17})
with initial trace $(\mathcal{S}(u),\xm_u).$ Then $u\geq u_{\xm_u}.$
\end{prop}
\Proof We assume $\mathcal{S}(u)\neq\BBR^n$ otherwise the result is proved. Let $G$ and $E$ be open bounded domains such that $G\subset\subset E\subset\subset \mathcal{R}(u).$ Let $0<\xd=\inf\{|x-y|:\;x\in G,\;y\in E^c\}.$ Choose $R>0$ such that $E\subset\subset B_R.$ Let $\{t_j\}_{j=1}^\infty$ be a decreasing sequence converging to $0$. We denote by $u_j$ the weak solution of the problem
$$\BA {lll}
\prt_tv-\xD v+Vv=0\qquad &\mathrm{in }\;\; B_R\times(t_j,T]\\
\phantom{--,,Vv}
v(x,t)=0\qquad&\mathrm{on }\;\; \partial B_R\times(t_j,T]\\
\phantom{--,,Vv}
v(.,t_j)=\chi_{G}u(.,t_j)\qquad &\mathrm{in }\;\;B_R\times\{t_j\},
\EA
$$
where $\chi$ is the characteristic function on $G.$
Let $v_j^R$, be the solution
$$\BA {lll}
\prt_tv-\xD v=0\qquad &\mathrm{in }\;\; \mathbb{R}^n\times(t_j,\infty]\\ \phantom{,-}
v(.,t_j)=\chi_{G}u(.,t_j)\qquad &\mathrm{in }\;\;\mathbb{R}^n\times\{t_j\}.
\EA
$$
Then by maximum principle we have $u_j^R\leq u$ and $u_j^R\leq v_j$  in $B_R\times(t_j,T]$, for any $j\in\mathbb{N}$. By standard parabolic estimates,  we may assume that the sequence $u_j^R$ converges locally uniformly in $Q_T^{B_R}$ to a function $u^R\leq u$. Moreover, since $\chi_{G}\xm_u(.,t_j)\rightharpoonup\chi_{G}\xm_u$ in the weak* topology, we derive from the  representation formula that $v_j\to  \BBH[\chi_{G}\xm_u]$. Furthermore $u^R\leq v$,  which implies $\chi_{(t_j,T)}u_j^R\to  u^R$ in $ L^1(Q_T^{B_R}).$ There also holds
$$\int_{t_j}^T\int_{B_R}u_j^RVdxdt=\int_{t_j}^T\int_{E}u_j^RVdxdt+\int_{t_j}^T\int_{B_R\setminus E}u_j^RVdxdt,$$
and, by the choice of $E$ and dominated convergence theorem,
$$\int_{t_j}^T\int_{E}u_j^RVdxdt\leq\int_{0}^T\int_{E}uVdxdt<\infty\Rightarrow
\lim_{j\to \infty}\int_{t_j}^T\int_{E}u_j^RVdxdt=\int_{0}^T\int_{E}u^RVdxdt.$$
Furthermore, for any $x\in B_R\setminus E$,
$$v_j(x,t)=\left(\frac{1}{4\pi(t-t_j)}\right)^\frac{n}{2}\int_{R^n}e^{-\frac{|x-y|^2}{4(t-t_j)}}\chi_{G}u(y,t_j)dy
\leq\left(\frac{1}{4\pi(t-t_j)}\right)^\frac{n}{2}e^{-\frac{\xd^2}{4(t-t_j)}}\int_{G}u(y,t_j)dy.$$
Next, since $V(x,t)\leq Ct^{-1}$ and $u_j^R\leq v_j,$ we obtain

\begin{equation}\label{Z1}
\lim_{j\to \infty}\int_{t_j}^T\int_{B_R\setminus E}u_j^RVdxdt=\int_{0}^T\int_{B_R\setminus E}u^RVdxdt,
\end{equation}
by using the previous estimate and the fact that $\chi_{G}\xm_u(x,t_j)\rightharpoonup\chi_{G}\xm_u$ in the weak* topology.
It follows $\chi_{(t_j,T)}Vu_j^R\to  Vu^R$ in $ L^1(Q_T^{B_R}).$
There holds also $u^R_G\leq u$; by the maximum principle,  the mapping $R\mapsto u^R_G$ is increasing and bounded from above by $u.$ In view of Lemma \ref{repr1},
$$\lim_{R\to \infty}u^R_G=u_G\leq u,$$
and $u_G$ is a positive weak solution of
$$\BA {lll}
\!\prt_tv-\xD v+Vv=0\qquad &\mathrm{in}\;\; Q_T\\ \phantom{\xD v+Vv}
v(.,0)=\chi_{G}\xm_u\qquad &\mathrm{in}\;\;\mathbb{R}^n.
\EA
$$
Consider an  increasing sequence $\{G_i\}_{i=1}^\infty$ of bounded open subsets, $G_i\subset \subset\CR(u)$, with the property that $\bigcup_{i=1}^\infty G_i=\mathcal{R}(u)$. In view of Lemma \ref{repr1} the sequence $\{u_i=u_{G_i}\}_{i=1}^\infty$ is increasing and converges to $\widetilde{u}\leq u.$ Also we have
$$u_i(x,t)+\int_0^t\int_{\mathbb{R}^n}H(t-s,x,y)V u_idyds=\int_{\mathbb{R}^n}H(x-y,t)d\xm_i,\;\forall (x,t)\in Q_T,$$
where $\xm_i=\chi_{G_i}\xm.$ Now since $\xm_i\rightharpoonup\xm_u,$ by the monotone convergence theorem we have
$$\tilde{u}(x,t)+\int_0^t\int_{\mathbb{R}^n}H(t-s,x,y)V \tilde{u}dyds=\int_{\mathbb{R}^n}H(x-y,t)d\xm_u,\;\forall (x,t)\in Q_T,$$
and $\tilde{u}\leq u.$ this implies $\tilde{u}=u_{\xm_u},$ which ends the proof.\hfill$\Box$\medskip

\noindent{\bf Remark}. Assumption (\ref{equ17}) is too strong and has only been used in  (\ref{Z1}). It could have been replaced by the following much weaker one: for any $R>0$ there exists a positive increasing function $\ge_R$ such that $\lim_{t\to 0}\ge(t)=0$ satisfying
\begin{equation}\label{Z2}
V(x,t)\leq e^{t^{-1}\ge_R(t)}\qquad\forall (x,t)\in Q_T^{B_R}.
\end{equation}

We end this section with a result which shows that the stability of the initial value problem with respect to convergence the initial data in the weak* topology implies that the initial of positive solution has no singular part.
\begin{theorem}
Assume $V$ satisfies, for some $\gt_0>0$,
\begin{equation}\label{Z3}
\lim_{\abs{E}\to 0}\dint_{E}H(x-y,t)V(x,t+\gt)dxdt=0,\qquad E\text{ Borel subset of }Q^{B_R}_T
\end{equation}
for any $R>0$, uniformly with respect to $y$ is a compact set and $\gt\in [0,\gt_0]$. If $u$ is a positive solution of
(\ref{equ17}), then $\mathcal{R}(u)=\mathbb{R}^n$
\end{theorem}
\Proof
We assume that $\mathcal{S}(u)\neq\emptyset$ and if $z\in \mathcal{S}(u)$
there holds
$$\dint_{\!\!Q^{B_r(z)}_T} Vudxdt=\infty\qquad\forall r>0.$$

In view of Proposition \ref{s(u)}, there exist two sequences $\{r_k\}$ and $\{t_j\}$ decreasing to $0$ such that
$$\lim_{j\to \infty}\int_{B_{r_k}(z)}u(x,t_j)dx=\infty\qquad\forall k\in\BBN.$$
For $k\in\BBN$ and $m>0$ fixed, there exists $j(k)$ such that
$$\int_{B_{r_k}(z)}u(x,t_j)dx\geq m\qquad\forall j\geq j(k)\in\BBN,$$
and there exists $\ell_{k}>0$ such that
$$\int_{B_{r_k}(z)}\min\{u(x,t_{j(k)}),\ell_{k}\}dx= m
$$
Furthermore $j(k)\to\infty$ when $k\to\infty$.
Let $R>\max\{r_j:j=1,2,...\}$ and  $u_k$ be the solution of
$$\BA {lll}
\prt_tv-\xD v+Vv=0\qquad &\mathrm{in }\;\; \BBR^n\times(t_{j(k)},T]\\
\phantom{--,,Vv}
v(.,t_j)=\chi_{B_{r_k}(z)}\min\{u(.,t_{j(k)}),\ell_{k}\}\qquad &\mathrm{in }\;\;\BBR^n\times\{t_{j(k)}\},
\EA
$$
Then $\chi_{B_{r_k}(z)}\min\{u(.,t_{j(k)}),\ell_{k}\}\to m\gd_z$ in the weak sense of measures.
By Proposition \ref{prop}  we obtain that $u\geq u_k$ on $B_R(z)\times (t_{j(k)},T]$. Applying Proposition \ref{uncondi+}, and the remark here after, we conclude that
$u_k(.,.+t_{j(k)})\to u_{m\gd_z}=mu_{\gd_z}$ in $L^1_{loc}(\overline Q_R^T)$ This implies $u\geq mu_{\gd_z}$, and as $m$ is arbitrary, $u=\infty$, contradiction.
\hfill$\Box$
\subsection{The sweeping method}
In this subsection we adapt to equation (\ref{equ17}) the sweeping method developed in \cite{VY} for constructing the boundary trace of solutions of stationnary Shr\"odinger equations. If $A\subset \BBR^n$ is a Borel set, we denote by
$${\mathfrak M}_{T\,+}(A)=\{\xm\in\mathfrak M_+(\BBR^n):\;  \gm(A^c)=0,\; \int_{A}e^{-\frac{|x|^2}{4T}}d\xm<\infty\}.$$
We recall that $\gm^*$ denotes the reduced measure associated to $\gm$.
\begin{prop}
Let $u\in W_{loc}(Q_T)$ be a positive solution of (\ref{equ17}) with singular set $\mathcal{S}(u)\varsubsetneq\mathbb{R}^n$. If $\xm\in{\mathfrak M}_{T\,+}(\mathbb{\mathcal{S}}(u))$, we set $v_\xm=\inf\{u,u_{\xm^*}\}.$ Then
$$\prt_t v_\xm-\xD v_\xm+Vv_\xm\geq 0\qquad \mathrm{in}\;\; Q_T,$$
and $v_\xm$ admits a boundary trace $\xg_u(\xm)\in\tilde{\mathfrak M}_+(\mathbb{\mathcal{S}}(u)).$ The mapping $\xm\mapsto\xg_u(\xm)$ is nondecreasing and $\xg_u(\xm)\leq\xm.$
\end{prop}
\Proof It is classical that $v_\gm:=\inf\{u,u_{\xm^*}\}$ is a supersolution of (\ref{equ17}) and $v_\gm\in \CE_\gn(Q_T)$ as it holds with $u_{\xm^*}$ .
The function
$$(x,t)\mapsto w(x,t)=\int_0^t\int_{\mathbb{R}^n} H(t-s,x,y)V(y,s)v_{\xm}(y,s)dyds$$
satisfies
$$\BA {ll}\prt_t w-\xD w-Vw= 0\quad \mathrm{in}\;\; Q_T\\\phantom{\prt_t w-\xD w}
w(.,0)=0\quad \mathrm{in}\;\; \BBR^n\times\{0\}.
\EA$$
Thus $v_\xm+w$ is a nonnegative supersolution of the heat equation in $Q_T$. It admits an initial trace in ${\mathfrak M}_{T\,+}(\mathbb{\mathcal{S}}(u))$ that we denote by $\xg_u(\xm).$ Clearly $\xg_u(\xm)\leq\xm^*\leq\xm$ since $v_\xm\leq u_{\xm^*}$ and $\xg_u(\xm)$ is nondecreasing with respect to $\xm$ as it is the case with $\xm\mapsto u_{\xm^*}$ is. Finally, since $v_\xm$ is a positive supersolution, it is larger that the solution of
\ref{sub1} where the initial data $\gm$ is replaced by $\xg_u(\xm)$, that is
$u_{\xg_u(\xm)}\leq v_\xm$.
\hfill$\Box$\\

The proofs of the next four propositions are mere adaptations to the parabolic case of similar results dealing with elliptic equations and proved in \cite{VY}; we omit them.
\begin{prop}
Let $$\xn_S(u):=\sup\{\xg_u(\xm):\;\xm\in{\mathfrak M}_{T\,+}(\mathbb{\mathcal{S}}(u))\}.$$
Then $\xn_S(u)$ is a Borel measure on $\mathbb{\mathcal{S}}(u).$
\end{prop}

\begin{defin}
The Borel measure $\xn(u)$ defined by
$$\xn(u)(A):=\xn_S(u)(A\cap\mathcal{S}(u))+\xm_u(A\cap\mathcal{R}(u)),\qquad\forall\;A\subset\mathbb{R}^n,\;A\;\mathrm{Borel},$$
is called the extended initial trace of $u,$ denoted by $tr_{\{t=0\}}^e(u).$
\end{defin}

\begin{prop}
If $A\subset\mathcal{S}(u)$ is a Borel set, then
$$\xn_S(A):=\sup\{\xg_u(\xm)(A):\;\xm\in{\mathfrak M}_{T\,+}(A)\}.$$
\end{prop}

\begin{prop}
There always holds $\xn(\mathcal{S}ing_V(\mathbb{R}^n))=0$, where $\mathcal{S}ing_V(\mathbb{R}^n)$ is defined
 in (\ref{Sing}).
\end{prop}

\begin{prop}
Assume $V$ satisfies condition (\ref{Z3}). If $u$ is a positive solution of (\ref{equ17}), then $tr^e_{\{t=0\}}(u)=\xm_u\in{\mathfrak M}_{T\,+}(\mathbb{R}^n).$
\end{prop}
\section{Appendix: the case of a bounded domain}
\setcounter{equation}{0}
\subsection{The subcritical case}
Let $\xO$ be a bounded domain with a $C^2$ boundary. We denote by  $\mathfrak M(\xO)$ the space of Radon measures in $\xO$, by $\mathfrak M_+(\xO)$ its positive cone and by $\mathfrak M_\gr(\xO)$ the space of Radon measures in $\xO$ which satisfy
\begin{equation}\label {W1}
\int_\Gw\gr d\abs{\gm}<\infty,
\end{equation}
for some weight function $\gr:\Gw\mapsto \BBR_+$. As an important particular case $\gr(x)=d^\ga(x)$, where $d(x)=\dist (x,\prt\Gw)$ and $\ga\geq 0$.
We consider the linear parabolic problem
\begin{equation}\label{W2}\BA {ll}
\prt_tu-\xD u+Vu=0,\quad &\mathrm{in }\;\; Q_T^\xO=\xO\times(0,T]\\ \phantom{\prt_tu-\xD u+V}
u=0\qquad&\mathrm{on}\;\; \partial_l Q_T^\xO=\partial\xO\times(0,T]\\\phantom{ \prt_tu+Vu}
u(.,0)=\xm\qquad &\mathrm{in}\;\;\xO.
\EA\end{equation}

\begin{defin}
We say that $\xm\in\mathfrak M_d(\xO)$ is a good measure if  the above problem has a weak solution $u$, i.e.
there exists a function $u\in L^1(Q_T^\xO),$ such that $Vu\in L^1_d(Q_T^\xO)$ which satisfies
\be
-\int_0^T\int_\xO u(\prt_t\xf+\xD\xf)dxdt+\int_0^T\int_\xO Vu\xf dxdt=\int_\xO\xf(x,0)d\xm,\label{weak}
\ee
$\forall \xf\in C^{1,1;1}(\overline{Q_T^\xO})$ which vanishes on $\partial_lQ_T^\xO$ and satisfies $\xf(x,T)=0$.
\end{defin}

\begin{defin}
Let $H^\xO(x,y,t)$ be the heat kernel in $\xO.$ Then we say that $\xm\in\mathfrak M_d(\xO)$ is a admissible measure if
$$||\BBH^\xO[|\xm|]||_{L^1(Q^\xO_T)}=\int_{Q_T^\xO}\left(\int_\xO H^\xO(x-y,t)d|\xm(y)|\right)V(x,t)\psi(x)dxdt<\infty.$$
\end{defin}

The next a proposition is direct consequence of \cite [Lemma 2.4] {MV}.

\begin{prop}\label{prop1}
Assume $\xm\in\mathfrak M_d(\xO)$ and  let $u$ be a weak solution of problem (\ref{W2}), then the following inequalities are valid\\
(i)
$$||u||_{L^1(Q_T^\xO)}+||Vu||_{L^1_\psi(Q_T^\xO)}\leq C(n,\xO)\int_\xO d d|\xm|,$$

\noindent (ii)
$$
-\int_0^T\int_\xO|u|(\prt_t\xf+\xD\xf)dxdt+\int_0^T\int_\xO |u|V\xf dxdt\leq\int_\xO\xf(x,0)d|\xm|,
$$
$\forall \xf\in C^{1,1;1}(\overline{Q_T^\xO})$, $\xf\geq 0$.\\
(iii)
$$
\xl_\xO\int_0^T\int_\xO (x)u^+dxdt+\int_0^T\int_\xO Vu^+\psi dxdt\leq\int_\xO\psi(x)d\xm^+.
,$$
where $\psi$ is the solution of
\begin{equation}\label{W3}\BA {ll}
-\xD \psi =1,\quad &\mathrm{in }\;\; \xO\\ \phantom{-\xD }
\psi=0\qquad&\mathrm{on}\;\; \partial \xO.
\EA\end{equation}
\end{prop}

\noindent\Proof For (ii), in \cite [Lemma 2.4, p 1456] {MV}, above from the relation (2.39), we can take $\tilde{\xz}=\xg(u)\xz$ for some $0\leq\xz\in C^{1,1;1}(\overline{Q_T^\xO})$, since $u=0$ on $\partial_lQ_T^\xO$.
For (iii) we consider (as in \cite[Remark 2.5]{MV}) $\xf(x,t)=t\psi(x)$.
The inequality holds by the same type of calculations as in \cite{VY}.\hfill$\Box$

\begin{prop}
The problem (\ref{W2}) admits at most one solution. Furthermore, if $\xm$ is admissible, then there exists a unique solution;  we denote it $u_\xm$. \label{subcrit}
\end{prop}

Similarly as Theorem \ref{uncondi} and Proposition \ref{uncondi}, we have the following stability results
\begin{prop}\label{prop}
(i) Assume that $V$ satisfies the stability condition
\be
\lim_{\abs{E}\to 0}\dint_E H^\xO(x,y,t)V(y,t)d(x) dydt=0,\quad\forall E\subset Q^\Gw_T,\, E\;\mathrm{ Borel}.\label{stability}
\ee
uniformly with respect to $y\in\Gw$. If $\{\xm_k\}$ is a bounded sequence in $\mathfrak M_d(\xO)$ converging to $\xm$ in the dual sense of $\mathfrak M_d(\xO)$, then $(u_{\xm_k},Vu_{\xm_k})$ converges to $(u_{\xm},Vu_{\xm})$ in $L^1(Q^\Gw_T)\times L_d^1(Q^\Gw_T)$.
\label{stab}
(ii) Furthermore if
\be
\lim_{\abs{E}\to 0}\dint_E H^\xO(x,y,t+\gt_n)V(y,t)d(x) dydt=0,\quad\forall E\subset Q^\Gw_T,\, E\;\mathrm{ Borel}.\label{stability+}
\ee
uniformly with respect to $y\in\Gw$ and $\gt_k\in [0,\gt_0]$ converges to $0$ and $\{\gm_k\}$ is in (i), then the solutions
$u_{\gt_k,\gm_k}$ of the shifted problem
\begin{equation}\label{equV++}\BA {ll}
\prt_tu-\Gd u+Vu=0\quad&\text{on }\;\;\Gw\times (\gt_k,T)\\
\phantom{\prt_tu-\Gd u+V}
u=0\quad&\text{on }\;\;\prt \Gw\times (\gt_k,T)
\\
\phantom{-\Gd u+..}
u(.,\gt_k)=\gm_k&\text{on }\;\;\Gw\times \{\gt_k\}
\EA\end{equation}
(extended by $0$ on $(0,\gt_k)$) converge to $u_{\gm}$ in $L^1_{d}(Q^\Gw_T)$, and  $\{Vu_{\gm_k}\}$ converges to $Vu_\gm$ in $L^1_{d}( Q^\Gw_T)$.
\end{prop}
\Proof
We can easily see that the measure $\xm_n$ is admissible and uniqueness holds; furthermore any admissible measure is a good measure  is a good measure as in Theorem \ref{fragma}, and
\[
\dint_{\!\!Q^\Gw_{T}}u_{\xm_n}dxds
+\dint_{\!\!Q^\Gw_{T}}u_{\xm_n}V\psi dxds
\leq C\int_{\xO}d\xm_n<C.
\]
The remaining of the proof is similar to the one of Theorem \ref{uncondi}.\hfill$\Box$

\subsection{The supercritical case}
\begin{lemma}
Let $\{\xm_n\}_{n=1}^\infty$ be an increasing sequence of good measures converging to some measure $\xm$ in the weak* topology, then $\xm$ is good.
\end{lemma}
\Proof
Let $u_{\xm_n}$ be the weak solution of (\ref{W2}) with initial data $\xm_n.$ Then
by  Proposition \ref{prop} -(iii), $\{u_{\xm_n}\}$ is an increasing sequence. By \ref{prop} -(i) the sequence $\{u_{\xm_n}\}$ is  bounded in $L^1(Q_T^\xO)$. Thus $u_{\xm_n}\to  u\in L^1(Q_T^\xO).$ Also by (iii) of
Proposition \ref{prop}, we have that $Vu_{\xm_n}\to  Vu$ in $L^1_\psi(Q_T^\xO).$ Thus we can easily prove that $u$ is a weak solution of (\ref{W2}) with $\xm$ as initial data.\hfill$\Box$\medskip

Let
\begin{equation}\label{ZOM}
Z_V^\xO=\{x\in\xO:\int_{Q_T^\xO}H^\xO(t,x,y)V(y)\psi(y)dy=\infty\}.\end{equation}
We note that, since $H^\xO(t,x,y)\leq H(x-y,t)$ for any bounded $\xO$ with smooth boundary, it holds
$Z_V^\xO\subset Z_V.$
By the same arguments as in \cite{VY} we can prove the following results
\begin{prop}
Let $\xm$ be an admissible positive measure. Then $\xm(Z_V^\xO)=0$
\end{prop}
\begin{prop}
Let $\xm\in\mathfrak M_{d\,+}(\xO)$ such that $\xm(Z_V^\xO)=0,$ then $\xm$ is good.\label{goodbound}
\end{prop}
\begin{prop}
Let $\xm\in\mathfrak M_{d\,+}(\xO)$ be a good measure. Then the following assertions are equivalent:\\
(i) $\xm(Z_V^\xO)=0.$\\
(ii) There exists an increasing sequence of admissible measures $\{\xm_n\}$ which converges to $\xm$ in the weak*-topology
\end{prop}

\end{document}